\documentclass[11pt]{amsart}
\usepackage{graphicx,cancel,xcolor}
\usepackage{amssymb,amsfonts,amsmath}
\usepackage{graphics,epsfig, verbatim}
\newcommand{\ds}{\displaystyle }

\newcommand{\BI}{\boldsymbol{\mathcal{I}}}

\newcommand{\BIT}{\boldsymbol{\mathcal{I}^n_{tot}}}

\newcommand{\BITp}{\boldsymbol{\mathcal{I}^{n+1}_{tot}}}

\begin{document}
\title[Maze solving ability of slime mold]{{On modeling Maze solving ability of slime mold via a hyperbolic model of chemotaxis}}

\author{G. Bretti} 
 \author{R. Natalini} 
\thanks{Istituto per le Applicazioni del Calcolo ``M.Picone'', Rome, Italy, g.bretti@iac.cnr.it, roberto.natalini@cnr.it}


\begin{abstract}
Many studies have shown that Physarum polycephalum slime mold is able to find the shortest path in a maze. In this paper we study this behavior in a network, using a hyperbolic model of chemotaxis. Suitable transmission and boundary conditions at each node are considered to mimic the behavior of such an organism in the feeding process.   Several numerical tests are presented for special network geometries to show the qualitative agreement between our model and the observed behavior of the mold.

\end{abstract}

\maketitle

 \\

\section{Introduction}

In recent years, the behavior of Physarum polycephalum, a ``many-headed'' slime mold with multiple nuclei, has been studied in depth and considered as a model organism representing amoeboid movement and cell motility. 
 The body of the plasmodium, representing the main vegetative phase of Physarum polycephalum,  contains a network
of tubes (pseudopodia), by means of which nutrients and chemical signals circulate throughout
the organism  \cite{N2000}; the flow in the tube is bidirectional, as it can be observed to switch back and forth.
The movement is characterized by cytoplasmic streaming, driven by rhythmic contractions of
organism, that sustains and reorganizes tubes: a larger flow leads to a wider tube. 

An example of the behavior of starved plasmodium searching for food has been shown in an  experiment performed by Nakagaki, Yamada, and T{\'o}th in \cite{N2000}. First, they built a maze on an agar surface. Therefore, they put over the agar surface some parts of the plasmodium, which started to spread and coalesced to form a single organism that filled the maze.
After that,  they placed agar blocks containing oatmeal at the exits and, after some hours,  they observed the dead ends of the plasmodium shrank (the tubes with a smaller flux ended to disappear) and the colony, after exploring all possible connections, produced the formation of a single thick pseudopodium spanning the minimum distance between the food sources. See for instance some popular videos freely available on the web, as for instance this one \cite{youtube}.\\
Then, to summarize, the main features in the evolution of the slime mold system are the following two steps: dead end cutting, and the selection of the solution path among the competitive paths.

This behavior can be applied to both path-finding in a maze and path selection in a
transport network. The ability of solving this kind of problems of an unicellular organism has been deeply investigated using various mathematical models. Some kinetic models were proposed and studied in \cite{Mi2007,Mi2008,KTN2007}. The hydrodynamic properties of  the transportation of flux by the tubes of plasmodium and the mechanisms under the rhythmic oscillation of Physarum were investigated in \cite{KTN2005},  and biologically inspired models for adaptive network construction were developed \cite{T2010} and \cite{W2011}. A mathematical model based on ordinary differential equations to describe the dynamics of Physarum was introduced in \cite{KTN2007}. It is based on a regulation mechanism, which balances the system between the thickness of a tube and the flux through it. In that paper, it was shown numerically that, in the asymptotic steady state of the model, the solution converges to the minimum-length solution between a source and a sink on any input graph. To show this property rigorously, it is necessary to prove that the globally
asymptotically stable equilibrium point corresponds to the shortest path connecting two exits of a network and the system has neither the oscillating nor chaotic solution. In \cite{Mi2007} there was proof of convergence for two parallel links, while the convergence of this model for any network to the shortest path connecting the source and the sink was analytically proved in \cite{bonifaci}. \\
Another approach is represented by modeling by using partial differential equations. In fact, the movement of individuals under the effect of a chemical stimulus of chemoattractant  has been widely studied in mathematics in the last two decades, see \cite{H, M, Pe}, and various models involving partial differential equations have been proposed. 
This process is known as chemotaxis: the organism will move with higher probability towards food sources - positive chemotaxis and, while eating, it will produce a chemical substance, usually a chemokine. The opposite situation
is that it moves away (negative chemotaxis), if it experience the presence of a harmful substance. 
The basic unknowns in these chemotactic models are the density of individuals and the concentrations of some chemical attractants. One of the most considered models of this type is the  Patlak-Keller-Segel system \cite{KS},
 where the evolution of the density of cells is described by a parabolic equation, and the concentration of  a chemoattractant is generally given by a parabolic or elliptic equation, depending on the different regimes to be described and on authors' choices. 
The analytical study of the KS model on networks was presented in the recent paper \cite{CC}, where the existence of a time global and spatially continuous solution for both the doubly parabolic and the parabolic-elliptic systems was proved.\\
By contrast, models based on hyperbolic/kinetic equations  for the evolution of the density of individuals, are characterized by a finite speed of propagation \cite{GambaEtAl03,Pe, FilbetLaurencotPerthame, DolakHillen03, gumanari}.

In this paper, we deal with a hyperbolic-parabolic model, originally considered in \cite{Segel77}, and later reconsidered in \cite{GreenbergAlt87}, which arises as a simple model for chemotaxis on a line:
\begin{equation}\label{hyper-gen}
 \left\{\begin{array}{ll}
  u_t +v_x =0,\\ 
  v_t +\lambda^2 u_x 
   =\chi \phi_x\,u-v,\\ 
   \phi_{t}-D\, \phi_{xx}=au-b\phi.
 \end{array}\right.
\end{equation}
The function $u$ is the density of individuals in the considered medium, $v$ is their averaged flux and $\phi$ denotes the density of chemoattractant. The individuals move at a constant speed $\lambda \geq 0$, changing their direction along the axis during the time. The positive constant $D$ is the diffusion coefficient of the chemoattractant; the positive coefficients $a$ and $b$, are respectively its production and degradation rates, and $\chi$ is the chemotactic sensitivity. Let us underline that the flux $v$ in model (\ref{hyper-gen}) corresponds to $v= -\lambda^2 u_x + \chi \phi_x u$ for the parabolic KS system. This system was analytically studied on the whole line and on bounded intervals in \cite{gumanari}, while an effective numerical approximation, the Asymptotic High Order (AHO) scheme, was proposed in \cite{NaRi}, see also \cite{Gosse10, Gosse11}.
The AHO scheme was introduced in order to balance correctly the source term with the differential part and avoid  an incorrect approximation of the density flux at equilibrium. These schemes are based on standard finite differences  methods, modified by a suitable treatment of the source terms, and they take into account, using a Taylor expansion of the solution, of the behavior of the solutions near non constant stationary states.

Model (\ref{hyper-gen}) can also be considered on a  network, where solutions  on each arc are coupled through transmission conditions ensuring the total continuity of the fluxes at each node, while the densities can have jumps. The first analytical study of this model on a network  was given in \cite{GuNa}, where a global existence result of solution for suitably small initial data was established. Concerning the numerical approximation, in \cite{BNR14},  we extended the AHO scheme of \cite{NaRi} to the case of networks. In this case, a particular attention was paid to the proper setting of conditions at internal and external nodes in order to guarantee the conservation of the total mass.

Here, inspired by the study on a doubly parabolic chemotaxis system presented in \cite{BGKS,CC}, where the one-dimensional KS model was extended to networks, we present a  numerical study of amoeboid movements, using the hyperbolic chemotaxis model \eqref{hyper-gen} on networks connecting two or more exits. The main difference with  previous works is not only the hyperbolicity of the equations describing the cell movement, but the transmission conditions, which in our case, unlike what was assumed in \cite{BGKS,CC},  are designed to allow for discontinuous densities at each node. Actually, this condition is believed to be more appropriate when considering a flow of individuals on a network, see for instance \cite{pic} and references therein and also the recent discussion in \cite{BKP}.

As in \cite{BNR14}, the numerical approximation of the hyperbolic part of the system is based on the AHO scheme, while the parabolic part is approximated with the Crank-Nicolson scheme. However, since here we include the modeling of the inflow of the mass and the introduction of food at the network exits, we assume inflow boundary conditions at external nodes both for the density of cells and the chemoattractant. Therefore we need to deal with more general Neumann boundary conditions with respect to  \cite{BNR14}, both for the hyperbolic and the parabolic part of the system.

After proposing a modified scheme, this scheme will be used to obtain numerical solutions for networks with special geometries, motivated by experimental observations.  Let us notice that the properties of dead-end cutting and shortest path selection in our model are represented by different mass distribution on the edges of the assigned network described as an oriented graph. \\

The paper is organized as follows.
In Section \ref{sec2} we describe the main setting of problem \eqref{hyper-gen}. Section \ref{sec3} is devoted to the description of the numerical approximation of the problem based on the AHO scheme of second order with a suitable discretization of the transmission and  boundary conditions for the case of non-null Neumann boundary conditions. In particular, a strong emphasis will be given to formulate the correct boundary conditions for the AHO schemes for the inflow and outflow cases.
Finally, in Section \ref{sec4}, we report some numerical experiments for  the distribution of densities within networks of different topologies, having in mind the laboratory experiments made with Physarum, in order to test the correspondence between our simulations and the real behavior of such an organism.
In this framework, a comparison between the results obtained with the two models, the one in \cite{BGKS} and the present one, is also provided, showing substantially the same asymptotic behavior of the solutions, even if the transmission conditions at each node are not the same, and so the transitory profiles.\\

\section{Analytical framework}\label{sec2}
First, let us represent a network as a connected graph. 
We define a connected  graph $G=(\mathcal{N},\mathcal{A})$ as formed by two finite sets, a set of $P$ nodes (or  vertices) $\mathcal{N}$  and  a set of $N$ arcs (or edges) $\mathcal{A}$,  such that an arc connects a pair of nodes. Since arcs are bidirectional the graph is non-oriented, but we need to fix an artificial orientation in order to fix a sign to the velocities. 
 The network is therefore  composed of "oriented" arcs and there are  two different types of intervals at a node  $p \in \mathcal{N}$ : incoming ones -- the set of these intervals is denoted by $I_{p}$ -- and outgoing ones  -- whose set  is denoted by $O_{p}$. For example, on the network depicted in Figure \ref{fig:1}, $1 \in I$ and $2,3 \in O$.
 We will also denote in the following by $I_{out}$ and $O_{out}$ the set of the arcs incoming or outgoing from the outer boundaries.
 The $N$ arcs of the network are parametrized as intervals $[0,L_i]$, $i=1,\ldots,N$, and for an incoming arc, $L_{i}$ is the abscissa of the node, whereas it is $0$ for an outgoing arc.\\

We consider system (\ref{hyper-gen}) on each arc and  rewrite it    
in diagonal variables for its hyperbolic part by setting
\begin{equation}
\label{cvar}
  u^{\pm}=\frac{1}{2}\left(u\pm \frac{v}{\lambda}\right).
  \end{equation} 
Here $\ds u^{+}$ and  $\ds u^{-}$ are the Riemann invariants of the system and $\ds u^{+}$ (resp. $\ds u^{-}$)  denotes the density of cells following the orientation of the arc (resp. the density of cells going in the opposite direction). This transformation is inverted by $ u=u^++u^-$  and  $v=\lambda(u^+-u^-)$, and yields:
\begin{equation}
\label{hyper-gen-diag}
\left\{
\begin{aligned}
u^+_{t}+\lambda u^+_{x}&=\frac{1}{2\lambda}\left((\chi\phi_{x}-\lambda)u^++(\chi\phi_{x}+\lambda)u^-\right),\\
u^-_{t}-\lambda u^-_{x}&=-\frac{1}{2\lambda}\left((\chi\phi_{x}-\lambda)u^++(\chi\phi_{x}+\lambda)u^- \right),\\
\phi_{t}-D\phi_{xx}&= a(u^++u^-)-b\phi.
\end{aligned}
\right.
\end{equation}
We complement this system by  initial conditions at $t=0$ on each arc
$$u^+(x,0) =u^+_{0}(x), \quad u^-(x,0) =u^-_{0}(x), \quad \phi(x,0)=\phi_{0}(x), \textrm{ for } x \in   [0,L], $$
with $u^+_{0}, \quad u^-_{0}, \quad \phi_{0}$ some  $C^3$ functions.
 
Up to now, we omitted the indexes related to the arc number since no confusion was possible. From now on, however, we need to distinguish the quantities on different arcs and we denote by $u_{i}^{\pm}$, $u_{i}$ , $v_{i}$,  $\phi_{i}$ etc., the values of the corresponding variables on the $i$-th arc.

On the outer boundaries, we consider the more general boundary conditions:
\begin{equation}\label{boundary}
\left\{
\begin{aligned}
u_{i}^+(0,t)=g_i(t, u_{i}^-(0,t)) ,&  \textrm{ if } i\in I_{out}, \\ 
u_{i}^-(L_{i},t)=g_i(t, u_{i}^+(L_i,t)), & \textrm{ if } i\in O_{out},
\end{aligned}
\right.
\end{equation}
for some nonlinear smooth functions $g_i(t, u)$,  which correspond to general boundary conditions on the flux function $v$:
 \begin{equation}\label{general_bound_v}
\left\{
\begin{aligned}
v_i(0,t) = \psi_i(t, u(0,t)), \textrm{ if } i \in I_{out},\\ 
v_i(L_i,t) = \psi_i(t, u(L_i,t)), \textrm{ if } i \in O_{out},\end{aligned}
\right.
\end{equation}
where $\psi_i(t,u)$ are the same functions after the change of variables.
We can also assume non-vanishing conditions for the flux of chemoattractant, to reproduce the flow of nutrient substances at the outer boundaries
\begin{equation}\label{general_bound_phi}
\partial_{x}\phi_{i}(.,t)= \bar\phi_{i}(t),
\end{equation}
 where $\bar\phi_{i}$ are given functions.

\subsection{Transmission conditions}

Let us describe how to define the conditions at each node; this is an important point, since the behavior of the solution will be very different according to the conditions we choose. Moreover, let us recall that the coupling between the densities on the arcs are obtained through these conditions. 
At   node $p \in \mathcal{N}$, we have to give values to the components such that the corresponding characteristics are going out of the node. Therefore, we consider the following transmission conditions at each node:
\begin{equation}\label{transmission}
\left\{
\begin{aligned}
u_{i}^-(L_{i},t)=\sum_{j \in I_{p}} \xi_{i,j} u_{j}^+(L_{j},t)+\sum_{j \in O_{p}} \xi_{i,j} u_{j}^-(0,t) ,&  \textrm{ if } i\in I_{p}, \\
u_{i}^+(0,t)=\sum_{j \in I_{p}} \xi_{i,j} u_{j}^+(L_{j},t)+\sum_{j \in O_{p}} \xi_{i,j} u_{j}^-(0,t), & \textrm{ if } i\in O_{p},
\end{aligned}
\right.
\end{equation}
where the constant $\xi_{i,j} \in [0,1]$ are the transmission coefficients: they represent the probability that a cell at a node decides to move from the $i-$th to the $j-$th arc of the network, also including the turnabout on the same arc. 
Let us notice that the condition differs when the arc is an incoming or an outgoing arc. Indeed, for an incoming (resp. outgoing) arc, the value of the function $u_{i}^+$ (resp. $u_{i}^-$)  at the node is obtained through the system and we need only to define $u_{i}^-$ (resp. $u_{i}^+$) at the boundary. 

 These  transmission conditions ensure the continuity of fluxes at each node, but not the continuity of the densities; however, the continuity of the fluxes is enough since it implies that we cannot loose nor gain any cells during the passage through a node. This is obtained using a condition mixing the transmission coefficients  $\xi_{i,j}$ and the velocities of the  arcs connected at each node $p$.
Fixing a node and denoting the velocities of the  arcs by  $\lambda_{i} \ge 0, \,  i \in I_p \cup O_{p}$, in order to have the flux conservation at  node $p$, which is given by: 
\begin{equation}\label{flux}
 \sum_{i \in I_{p} } \lambda_{i} (u_{i}^+(L_{i},t)-u_{i}^-(L_{i},t))=\sum_{i \in O_{p} } \lambda_{i} (u_{i}^+(0,t)-u_{i}^-(0,t)), 
\end{equation}
it is enough to impose the following condition:  
\begin{equation}\label{condition_lambda}
\sum_{i \in I_p \cup O_{p}} \lambda_{i}  \xi_{i,j} =\lambda_{j},  \, j \in I_p \cup O_{p}.
\end{equation}
 
The coefficients $\xi_{i,j}$ belong to $[0,1]$, and at every node $p \in \mathcal{N}$, we have:
\begin{equation}\label{dissip-coeff}
 \sum_{j  \in I_p \cup O_{p}}\xi_{i,j} = 1\textrm{ for all } i  \in I_p \cup O_{p}.
\end{equation}
In the case of no-flux at the outer boundaries, condition \eqref{condition_lambda} ensures that  the global mass $\mu(t)$  of the system is conserved along the time, namely:
\begin{equation}\label{masscons}
\mu(t) = \sum_{i=1}^N \int_{0}^{L_{i}} u_{i}(x,t) dx = \mu_0 := \sum_{i=1}^N \int_{0}^{L_{i}} u_{i}(x,0) dx, \, \textrm{ for all } t>0.
\end{equation}
In the case of non-vanishing Neumann conditions at the outer boundaries the conservation of the global mass of the system obviously does not hold. If we consider the case treated in the numerical tests in Section \ref{sec4}, we have two arcs connected to the outer boundaries (a source and a sink), say arcs $l$ and $m$. In such a case we have that the derivative in time of the total mass $\mu$ is given by:
\begin{equation}\label{massnotcons}
\mu'(t)  = v_l(0,t) - v_m(L_m,t), \, \textrm{ for all } t>0.
\end{equation}

Now let us consider the transmission conditions for $\phi$ in  system \eqref{hyper-gen}. We complement conditions \eqref{boundary} and \eqref{transmission} with a transmission condition  for $\phi$. 
As previously, we do not impose the continuity of the  density of chemoattractant  $\phi$, but only  the continuity of the flux at  node $p \in \mathcal{N}$. Therefore, we use  the Kedem-Katchalsky permeability condition \cite{KK}, which has been first proposed in the case of flux through a membrane. For some positive coefficients $\kappa_{i,j}$, we impose at each node
\begin{equation}
D_{i}\partial_{n}\phi_{i}=\sum_{j \in I_{p} {\cup} O_{p}} \kappa_{i,j} (\phi_{j}-\phi_{i}), \  i\in I_{p}\cup O_{p},
\label{transmission-phi}
\end{equation}
that, under the condition $\kappa_{i,j}= \kappa_{j,i}, i,j=1,\ldots,N$ the  conservation of the fluxes at  node $p$ holds, that is to say:
\[\ds \sum_{i \in I_{p} {\cup} O_{p}} D_{i}\partial_{n}\phi_{i}=0. \]
Note that setting $\ds \kappa_{i,i}=0, \, i=1,\ldots,N$, does not change condition \eqref{transmission-phi}, and the positivity of the transmission coefficients  $\kappa_{i,j}$, guarantees the energy dissipation for the equation for $\phi$ in  \eqref{hyper-gen}, when the term in $u$ is absent.

\section{A numerical approximation for system \eqref{hyper-gen}}\label{sec3}
Here we recall the numerical scheme presented in \cite{BNR14} and adapt this  scheme  to the case of non-vanishing Neumann boundary conditions.

\subsection{Approximation of the hyperbolic equations of the system in the case of a network.}
 Let us consider a network as previously defined. Each arc $a_{i}\in \mathcal{A}, \, 1\leq i \leq N$, is parametrized as an interval $a_{i}=[0,L_{i}]$ and is discretized with a space step $h_{i}$ and discretization  points $x^{j}_{i}$ for $j=0,\ldots, M_i+1$. We still denote by $k$ the time step, which is the same for all the arcs of the network.
In this subsection, we denote by $w^{n,j}_{i}$ the discretization on the grid at time $t_{n}$ and at point $x^{j}_{i}$  of a function $w_i, \ i=1,\ldots,N$ on the $i$-th arc for $j=0,\ldots, M_i+1$ and $n\ge 0$.
Here we describe the discretization of system \eqref{hyper-gen} on each arc, denoting by $f= \chi \phi_x\,u$ and omitting the parabolic equation for $\phi$. 
 Since we also work with Neumann boundary conditions for the $\phi$ function, the function $f$ will satisfy the following conditions on the boundary~:
$\ds f(0,t)=f(L,t)=0$. We therefore consider the following system  
\begin{equation}\label{hypergen2}
 \left\{\begin{array}{ll}
  u_t +v_x =0,\\ 
  v_t +\lambda^2 u_x 
   =f-v,
 \end{array}\right.
\end{equation}
and rewrite it in a   diagonal  form, using the  usual change of variables  \eqref{cvar},
\begin{equation}
\label{systeme}
\left\{
\begin{aligned}
u^-_{t}-\lambda u^-_{x}&=\frac{1}{2}(u^+-u^-)-\frac{1}{2\lambda}f, \\
u^+_{t}+\lambda u^+_{x}&=\frac{1}{2}(u^--u^+)+\frac{1}{2\lambda}f.
\end{aligned}
\right.
\end{equation}

To have a reliable scheme, with a correct resolution of fluxes at equilibrium, we have to deal with  Asymptotically High Order schemes on each arc $i$, see \cite{NaRi} for the general conditions on such coefficients in order to have a consistent, monotone scheme, that is of higher order on the stationary solutions.\\
In particular, the second order scheme in diagonal variables reads as:
\begin{eqnarray}
u_{-,i}^{n+1,j}&= \left(1-\lambda_{i}\frac{k}{h_{i}}-\frac{k}{4}\right)u_{-,i}^{n,j} +  \left(\frac{k \lambda_{i}}{h_{i}} - \frac{k}{4}\right) u_{-,i}^{n,j+1} + \frac{k}{4} (u_{+,i}^{n,j}  + u_{+,i}^{n,j+1})\nonumber\\
& - \frac{k}{4\lambda_{i}}\bigl(f_{i}^{n,j+1} + f_{i}^{n,j}\bigr), \quad j=0,\ldots,M_i, \label{scheme-aho_ord2:1} \\
u_{+,i}^{n+1,j}&= \left(1-\lambda_{i}\frac{k}{h_{i}}-\frac{k}{4}\right)u_{+,i}^{n,j}+  \left(\frac{k \lambda_{i}}{h_{i}} - \frac{k}{4}\right) u_{+,i}^{n,j-1} +\frac{k}{4} (u_{-,i}^{n,j} + u_{-,i}^{n,j-1})\nonumber\\
& + \frac{k}{4\lambda_{i}}\bigl(f_{i}^{n,j-1} + f_{i}^{n,j}\bigr), \quad j=1,\ldots,M_i+1. \label{scheme-aho_ord2:2}
\end{eqnarray}
Note that we have a second--order  AHO scheme on every stationary solutions, which is enough to balance the flux of the system at equilibrium.
It is possible to see that component by component  monotonicity conditions (neglecting the source term $f$)  are satisfied if 
\begin{equation}\label{cfl} h\leq 4 \lambda \text{ and } \ds k \leq \frac{4h}{h+4\lambda},\end{equation}
see  \cite{NaRi} for more details. For simplicity here we only deal with second order schemes, but third order schemes could also be considered (these schemes require a  fourth--order  AHO scheme for the parabolic equation with a five-points discretization for $\phi_{x}$). 
Finally, we rewrite the second order scheme \eqref{scheme-aho_ord2:1}-\eqref{scheme-aho_ord2:2} in the  variables $u$ and $v$  as~:
  \begin{equation}
\label{scheme-uv-interval}
\begin{aligned}
&u^{n+1,j}_{i}=u^{n,j}_{i}-\frac{k}{2h_{i}} \left(v^{n,j+1}_{i}-v^{n,j-1}_{i}\right)+ \frac{\lambda_{i}k}{2h_{i}}(u^{n,j+1}_{i}-2u^{n,j}_{i}+u^{n,j-1}_{i})\\ 
& +\frac{k}{4\lambda_{i}} (-v^{n,j-1}_{i} + v^{n,j+1}_{i}) +\frac{k}{4\lambda_{i}} (f^{n,j-1}_{i} - f^{n,j+1}_{i}),\\
&  v^{n+1,j}_{i}=v^{n,j}_{i} -\frac{\lambda^2_{i} k}{2h_{i}} \left(u^{n,j+1}_{i}-u^{n,j-1}_{i}\right)+ \frac{\lambda_{i}k}{2h_{i}}(v^{n,j+1}_{i}-2v^{n,j}_{i}+v^{n,j-1}_{i})\\
& -\frac{k}{4} (v^{n,j-1}_{i} + 2v^{n,j}_{i} + v^{n,j+1}_{i})
+\frac{k}{4\lambda_{i}} (f^{n,j-1}_{i} + 2f^{n,j}_{i} + f^{n,j+1}_{i}).
\end{aligned}
\end{equation}

\subsubsection{Boundary conditions}
Now, we define the numerical boundary conditions associated to this scheme.
As shown in \cite{BNR14},  to implement this scheme on each interval we need two boundary conditions at the endpoints linked to external nodes and two transmission conditions  at the endpoints linked to internal nodes. 
Considering an arc and its initial and end nodes, there are two possibilities: either they are external nodes, namely nodes from the outer boundaries linked to only one arc, or they are internal nodes connecting several arcs together. The boundary and transmission conditions will therefore depend on this feature.  Below, we will impose the boundary conditions for $v$:  \eqref{bound1} or \eqref{bound_notnull} and \eqref{bound3}   at outer nodes,  and two transmission conditions \eqref{bound2}--\eqref{bound4}  at inner nodes.

The first type of boundary conditions for will come from no-flux condition  at outer nodes~:
\begin{equation}\label{no-flux}
\left\{
\begin{aligned}
v_{i}^{n+1,0}=0, \textrm{ if } i\in I_{out}, \\
v_{i}^{n+1,M_{i}+1}=0, \textrm{ if } i\in O_{out},
\end{aligned}\right.
\end{equation} 
where $I_{out}$  (resp.  $O_{out}$) means that the arc is incoming from (resp. outgoing to) the outer boundary. 
In the $u^{\pm}$-variables, these conditions become~:
\begin{equation}\label{bound1}
\left\{
\begin{aligned}
u_{+,i}^{n+1,0}=u_{-,i}^{n+1,0}, \textrm{ if } i\in I_{out}, \\
u_{+,i}^{n+1,M_{i}+1}=u_{-,i}^{n+1,M_{i}+1}, \textrm{ if } i\in O_{out}.
\end{aligned}\right.
\end{equation} 

The boundary condition at a node $p$ will come from a discretization of the transmission condition \eqref{transmission}, that is to say
 \begin{equation}\label{bound2}
\left\{
\begin{aligned}
u_{-,i}^{n,M_{i}+1}=\sum_{j \in I_{p}} \xi_{i,j} u_{+,j}^{n,M_{j}+1}+\sum_{j \in O_{p}} \xi_{i,j} u_{-,j}^{n,0} ,&  \textrm{ if } i\in I_{p}, \\
u_{+,i}^{n,0}=\sum_{j \in I_{p}} \xi_{i,j} u_{+,j}^{n,M_{j}+1}+\sum_{j \in O_{p}} \xi_{i,j} u_{-,j}^{n,0} , & \textrm{ if } i\in O_{p}.
\end{aligned}
\right.
\end{equation}
However, these relations link all the unknowns together and they cannot be used alone. An effective way to compute all these quantities will be presented after equation \eqref{bound4} below.
We still have two missing conditions per arc, which can be recovered by using the boundary conditions on $v$.

In particular, in the case of null Neumann boundary conditions (\ref{bound1}), from the exact mass conservation between two successive computational steps,
we obtain for the (\ref{scheme-aho_ord2:1})-(\ref{scheme-aho_ord2:2}) the following boundary conditions:
\begin{equation}
 \label{bound3}
 \left\{
\begin{aligned}
& u_{+,i}^{n+1,0}=u_{-,i}^{n+1,0}=(1-\lambda_{i}\frac{k}{h_{i}})u_{-,i}^{n,0}+\frac k 4 (\frac{\lambda_{i}}{h_{i}}-1) u_{-,i}^{n,1}
\\
& \qquad \qquad \qquad \qquad+\frac k 4 u_{+,i}^{n,1}-\frac{k}{4\lambda_{i}}f_{i}^{n,1},
   \textrm{ if } i \in I_{out},\\
& u_{+,i}^{n+1,M+1}=u_{-,i}^{n+1,M+1}=(1-\lambda_{i}\frac{k}{h_{i}})u_{+,i}^{n,M+1}+\frac k 4 (\frac{\lambda_{i}}{h_{i}}+ 1) u_{+,i}^{n,M}\\
&\qquad \qquad \qquad \qquad-\frac k 4 u_{-,i}^{n,M}-\frac{k}{4\lambda_{i}} f_{i}^{n,M_{i}}, 
\textrm{ if } i \in O_{out},
\end{aligned}
\right.
\end{equation}
where $I_{out}$ and $O_{out}$ have the same meaning as previously. 

Using the transmission conditions \eqref{bound2} for $\ds u_{-,i}^{n+1,M_{i}+1}$ if $i\in I_{p}$ and for $\ds u_{+,i}^{n+1,0}$ if $i\in O_{p}$, we obtain the following numerical boundary conditions, with $\delta_i=h_{i}\left(h_{i}+  \sum_{j \in I_{p}  {\cup} O_{p}} h_{j}\xi_{j,i}\right)^{-1}$:
\begin{equation}
 \label{bound4}
\begin{aligned}
& u_{+,i}^{n+1,M_{i}+1}= \delta_i \Bigl(u_{+,i}^{n,M_{i}+1}( 1-\frac{k}{2})
+u_{-,i}^{n,M_{i}+1}\bigl( 1-2 k\frac{\lambda_{i}}{h_{i}} + \frac{k}{2}\bigr)\\
&+k u_{+,i}^{n,M_{i}}\bigl(2\frac{\lambda_{i}}{h_{i}}-\frac{1}{2}\bigr)
-\frac{k}{2} u_{-,i}^{n,M_{i}} + \frac{k}{\lambda_{i}}( \frac{1}{2}f_{i}^{n,M_{i}+1} +\frac{1}{2} f_{i}^{n,M_{i}})\Bigr), \textrm{ if } i\in I_{p},
\end{aligned}
\end{equation}
\begin{equation*}
\begin{aligned}
& u_{-,i}^{n+1,0}= \delta_i \Bigl(u_{+,i}^{n,0}( 1-2k \frac{\lambda_{i}}{h_{i}}
+\frac{k}{2}) +u_{-,i}^{n,0}\bigl( 1- \frac{k}{2}\bigr)+  \frac{k}{2} u_{+,i}^{n,1}\\
&+k u_{-,i}^{n,1}\bigl(2\frac{\lambda_{i}}{h_{i}}-\frac{1}{2}\bigr)
-\frac{k}{\lambda_{i}}(\frac{1}{2}f_{i}^{n,1}+\frac{1}{2}f_{i}^{n,0}) 
\Bigr)
, \textrm{ if } i\in O_{p}.
\end{aligned}
\end{equation*}
Once these quantities are computed, we can use equations \eqref{bound2} at time $t_{n+1}$, to obtain $u_{-,i}^{n+1,M_{i}+1}$ if $i \in I_{p}$  and $u_{+,i}^{n+1,0}$ if $i \in O_{p}$.

As proved in \cite{BNR14}, the consistency of (\ref{bound4}) at a node holds under the condition on each arc 
\begin{equation}\label{CFL}
h_i = 2 k \lambda_i.
\end{equation}
Please notice that using this condition in combination with \eqref{cfl}, implies that $k\leq 2$. 
In conclusion, we have imposed four boundary conditions \eqref{bound1},  \eqref{bound2}, \eqref{bound3}, and \eqref{bound4} on each interval. Conditions \eqref{bound1} and \eqref{bound3} deal with the outer boundary, whereas conditions \eqref{bound2} and \eqref{bound4} deal with the node. 
Under these conditions, the total numerical mass is conserved at each step. 

Let us now consider the more general case of non vanishing Neumann boundary conditions, which is finally the main goal of this section. In particular, we consider, as in \cite{BGKS}, the inflow condition at the nodes $p_1$ and $p_2$, respectively: 
\begin{equation}\label{inflow_cond_v1} 
v_l(0,t) = \frac{2}{1+ u_l(0,t)}, \ l \in O_{p_1}
\end{equation}
\begin{equation}\label{inflow_cond_v2} 
\quad v_{m}(L_{m},t) = -\frac{2}{1+ u_{m}(L_{m},t)} \ m \in I_{p_2}.
\end{equation}
 The discretization of conditions (\ref{inflow_cond_v1})-(\ref{inflow_cond_v2}) is:
\begin{equation}\label{bound_notnull}
\left\{
\begin{aligned}
v_{l}^{n+1,0}&= \frac{2}{1 + u_{+,l}^{n+1,0}+u_{-,l}^{n+1,0}},  \ l \in O_{p_1}\\
v_{m}^{n+1,M_{m}+1}&=-\frac{2}{1 + u_{+,m}^{n+1,M_{m}+1}+u_{-,m}^{n+1,M_{m}+1}}, \ m \in I_{p_2}.
\end{aligned}\right.
\end{equation} 

Since in the general case of non-null Neumann conditions the conservation of the total mass of the system does not hold, we need to compute the numerical approximation scheme at the outer boundaries in such a case, taking into account the relation (\ref{massnotcons}).
 The discrete total mass is given by $\ds \BIT=\sum_{i=1}^N \BI_{i}^n$,  where the mass corresponding to the arc $i$ is defined as:
\begin{equation*}
\begin{aligned}
 \BI_{i}^n&=h_{i}\left(\frac{u_{i}^{n,0}}{2}+\sum_{j=1}^{M_{i}} u_{i}^{n,j}+\frac{u_{i}^{n,M_{i}+1}}{2}\right)\\
 &= h_{i}\left(\frac{u_{+,i}^{n,0}+u_{-,i}^{n,0}}{2}+\sum_{j=1}^{M_{i}} (u_{+,i}^{n,j}+u_{-,i}^{n,j})+\frac{u_{+,i}^{n,M_{i}+1}+u_{-,i}^{n,M_{i}+1}}{2}\right).
\end{aligned}
\end{equation*}
Computing $\ds \BITp-\BIT$, we find:
\begin{equation*}
\begin{split}
&\BITp-\BIT=\sum_{i=1}^N \frac{h_{i}k}{2} 
\Biggl( \frac{1}{k}(u_{+,i}^{n+1,0}-u_{+,i}^{n,0})+ \frac{1}{k}(u_{-,i}^{n+1,0}-u_{-,i}^{n,0})
+ (2\frac{\lambda_{i}}{h_{i}}-  \frac{1}{2} )u_{+,i}^{n,0}\\
&+ \frac{1}{2} u_{-,i}^{n,0}- \frac{1}{2} u_{+,i}^{n,1}-(2\frac{\lambda_{i}}{h_{i}}- \frac{1}{2})u_{-,i}^{n,1}
 +\frac{1}{\lambda_{i}}\left(\frac{1}{2}f_{i}^{n,1}+ \frac{1}{2}f_{i}^{n,0}\right) 
\Biggr)\\
&+\frac{h_{i}k}{2} \Biggl( \frac{1}{k}(u_{+,i}^{n+1,M_{i}+1}-u_{+,i}^{n,M_{i}+1})+\frac{1}{k}(u_{-,i}^{n+1,M_{i}+1}-u_{-,i}^{n,M_{i}+1})
  + \frac{1}{2} u_{+,i}^{n,M_{i}+1}\\
  &+(2\frac{\lambda_{i}}{h_{i}}- \frac{1}{2})u_{-,i}^{n,M_{i}+1} -(2\frac{\lambda_{i}}{h_{i}} - \frac{1}{2})u_{+,i}^{n,M_{i}} -\frac{1}{2} u_{-,i}^{n,M_{i}}\\
& - \frac{1}{\lambda_{i}}\left(\frac{1}{2} f_{i}^{n,M_{i}+1} + \frac{1}{2} f_{i}^{n,M_{i}}\right)
\Biggr).
\end{split}
\end{equation*}
We are going to impose boundary conditions such that the right-hand side in the previous difference is equal to the difference between the fluxes at the outer boundaries. In particular, in the case of the networks considered in Section \ref{sec4}, at the outer boundaries we have the inflow conditions \eqref{inflow_cond_v1} for an arc $l$ exiting from the source node and \eqref{inflow_cond_v2} for an arc $m$ entering a sink node. Then, we discretize the relation (\ref{massnotcons}), and, using the trapezoidal rule it can be written as: 
$$\BITp-\BIT= \frac{k}{2} (v_l^{n,0} - v_{m}^{n,M_{m}+1} + v_l^{n+1,0} - v_{m}^{n+1,M_{m}+1})$$
and then, splitting the terms, for the left boundary (on arc $l$) and right boundary (on arc $m$), we get, respectively,
\begin{equation}
\begin{split}\label{arcl}
&\frac{k}{2} \left(\frac{2}{1 + u_{+,l}^{n,0}+u_{-,l}^{n,0}}+ \frac{2}{1 + u_{+,l}^{n+1,0}+u_{-,l}^{n+1,0}}\right)\\
&=\frac{h_{1}k}{2} 
\Biggl( \frac{1}{k}(u_{+,l}^{n+1,0}-u_{+,l}^{n,0})+ \frac{1}{k}(u_{-,l}^{n+1,0}-u_{-,l}^{n,0})
+ (2\frac{\lambda_{l}}{h_{l}}-  \frac{1}{2} ) u_{+,l}^{n,0}\\
& + \frac{1}{2} u_{-,l}^{n,0}- \frac{1}{2} u_{+,l}^{n,1}-(2\frac{\lambda_{l}}{h_{l}}- \frac{1}{2})u_{-,l}^{n,1}
 +\frac{1}{\lambda_{l}}\left(\frac{1}{2}f_{l}^{n,1}+ \frac{1}{2}f_{l}^{n,0}\right)
\Biggr)
\end{split}
\end{equation}
and
\begin{equation}
\begin{split}\label{arcm}
& \frac{k}{2} \left(\frac{2}{1 + u_{+,m}^{n,M_{m}+1}+u_{-,m}^{n,M_{m}+1}} + \frac{2}{1 + u_{+,m}^{n+1,M_{m}+1}+u_{-,m}^{n+1,M_{m}+1}}\right)\\
&= \frac{h_{m}k}{2} \Biggl( \frac{1}{k}(u_{+,m}^{n+1,M_{m}+1}-u_{+,m}^{n,M_{m}+1})+\frac{1}{k}(u_{-,m}^{n+1,M_{m}+1}-u_{-,m}^{n,M_{m}+1})
  + \frac{1}{2} u_{+,m}^{n,M_{m}+1}\\
&+(2\frac{\lambda_{m}}{h_{m}}- \frac{1}{2})u_{-,m}^{n,M_{m}+1} -(2\frac{\lambda_{m}}{h_{m}} - \frac{1}{2})u_{+,m}^{n,M_{m}} -\frac{1}{2} u_{-,m}^{n,M_{m}}\\
& - \frac{1}{\lambda_{m}}\left(\frac{1}{2} f_{m}^{n,M_{m}+1} + \frac{1}{2} f_{m}^{n,M_{m}}\right)
\Biggr).
\end{split}
\end{equation}
Let us now focus on the condition for arc $l$.

Setting 
\begin{equation}\label{paramA}
\begin{split}
A^{n} &=- \frac{2}{1 + u_{+,l}^{n,0}+u_{-,l}^{n,0}}   + h_{l} \Biggl(- \frac{1}{k} (u_{+,l}^{n,0} +u_{-,l}^{n,0})
+ (2\frac{\lambda_{l}}{h_{l}}-  \frac{1}{2} ) u_{+,l}^{n,0}\\
& + \frac{1}{2} u_{-,l}^{n,0}- \frac{1}{2} u_{+,l}^{n,1}-(2\frac{\lambda_{l}}{h_{l}}- \frac{1}{2})u_{-,l}^{n,1}
 +\frac{1}{\lambda_{l}}\left(\frac{1}{2}f_{l}^{n,1}+ \frac{1}{2}f_{l}^{n,0}\right)
 \Biggr),
\end{split}
\end{equation}
and the other quantities
\begin{equation}\label{param1}
  \alpha_1 =\frac{h_l}{k}, \ \beta_1^n = \alpha_1 + A^{n}, \ \gamma_1^n = A^{n} -2,
\end{equation}
we can write the equation (\ref{arcl}) as
\begin{equation}\label{u1+_2}
  \alpha_1 (u_{+,l}^{n+1,0})^2 + (\beta_1^n + 2\alpha_1 u_{-,l}^{n+1,0}) u_{+,l}^{n+1,0} + \gamma_1^n + u_{-,l}^{n+1,0}(\alpha_1 (1+u_{-,l}^{n+1,0}) + A^n) =0.
\end{equation}
Then, choosing the positive root of equation (\ref{u1+_2}), we get the formula:
\begin{equation}\label{u1+}
\begin{split}
u_{+,l}^{n+1,0} &:=g_1^n(u_{-,l}^{n+1,0})=\frac{1}{2 \alpha_1} \Biggl[-\beta_1^n -2\alpha_1 u_{-,l}^{n+1,0}  \\
 &+\sqrt{\left(\beta_1^n + 2\alpha_1 u_{-,l}^{n+1,0}\right)^2 - 4 \alpha_1 \left(\gamma_1^n + u_{-,l}^{n+1,0} (A^{n} + \alpha_1) + \alpha_1 (u_{-,l}^{n+1,0})^2\right)}\Biggr].
\end{split}
\end{equation}
We compute (\ref{u1+}) using the value of $u_{-,l}^{n+1,0}$ obtained from the numerical scheme \eqref{scheme-aho_ord2:1} at the boundary $j=0$,
with
 $$
 f_l^{n,0}= \chi u_l^{n,0} \phi_{x,l}^{n,0}, 
$$
and
$$
 f_l^{n,1}=\chi u_l^{n,1} \phi_{x,l}^{n,1}= \chi u_l^{n,1} \frac{\phi_l^{n,2}-\phi_l^{n,0}}{2 h_l}.  
$$
Then, plugging the expression above into \eqref{scheme-aho_ord2:1} and denoting $u_{+,l}^{n,0}= g^{n-1}_1(u_{-,l}^{n,0})$, we have:
\begin{equation}\label{u1-_boundary:2}
\begin{split}
u_{-,l}^{n+1,0}&= u_{-,l}^{n,0}\left(1 - \frac{\lambda_l k}{h_l} - \frac{k}{4} - \frac{k}{4 \lambda_l} \chi \phi_{x,l}^{n,0}\right) + g^{n-1}_1(u_{-,l}^{n,0}) \left(\frac{k}{4} - \frac{k}{4 \lambda_l} \chi \phi_{x,l}^{n,0}\right)  \\
&+ u_{-,l}^{n,1} \left(\frac{\lambda_l k}{h_l}- \frac{k}{4} - \frac{k\chi}{4 \lambda_l} \frac{\phi_l^{n,2}-\phi_l^{n,0}}{2 h_l}\right) +  u_{+,l}^{n,1}\left(\frac{k}{4} - \frac{k \chi}{4\lambda_l} \frac{\phi_l^{n,2}-\phi_l^{n,0}}{2 h_l}\right) . 
\end{split}
\end{equation}
Let us now analyze the sign of the coefficients in (\ref{u1-_boundary:2}) in order to ensure the monotonicity of the scheme.
Under the condition (\ref{CFL}), we want to satisfy the following relations:

\begin{equation}\label{coeff_zw}
\left\{\begin{aligned}
&\frac{\lambda_l k}{h_l}- \frac{k}{4} - \frac{k\chi}{4 \lambda_l} \frac{\phi_l^{n,2}-\phi_l^{n,0}}{2 h_l} \ge 0,\\
&\frac{k}{4} - \frac{k \chi}{4\lambda_l} \frac{\phi_l^{n,2}-\phi_l^{n,0}}{2 h_l} \ge 0,
\end{aligned}\right.
\end{equation}
and we also have to guarantee the monotonicity with respect to $u_{-,l}^{n,0}$ of the first two terms on the right-hand side of equation \eqref{u1-_boundary:2}:
\begin{equation}\label{coeff_w0}
u_{-,l}^{n,0}\left(1 - \frac{\lambda_l k}{h_l} - \frac{k}{4} - \frac{k}{4 \lambda_l} \chi \phi^{n,0}_{x,l}\right) + g^{n-1}_1(u_{-,l}^{n,0}) \left(\frac{k}{4} - \frac{k}{4 \lambda_l} \chi \phi^{n,0}_{x,l}\right).
\end{equation}
For (\ref{coeff_zw}) we find the condition:
\begin{equation}\label{cond_coeff_zw}
 \phi_{x,l}^{n,1} \le \left(\frac{1}{k}-\frac{1}{2}\right)\frac{2\lambda_l}{\chi}, \ \textrm { with } k \le 1.
\end{equation}
It is not obvious to have that the condition (\ref{cond_coeff_zw}) is respected, since, in general, the gradient of $\phi$ can increase during chemotactic process if the mass grows, thus causing blow-up of solutions. Then at each time step we need to check if the condition is satisfied, in order to have a finite solution.
 Setting $u_{-,l}^{n,0}=u$ in (\ref{coeff_w0}) and passing to the derivative respect to $u$, we study the expression
\begin{equation}\label{der_coeff_w0}
\left(1 - \frac{\lambda_l k}{h_l} - \frac{k}{4} - \frac{k}{4 \lambda_l} \chi \phi_{x,l}^{n,0}\right) + (g_1^{n-1})'(u) \left(\frac{k}{4} - \frac{k}{4 \lambda_l} \chi \phi_{x,l}^{n,1}\right) \ge 0.
\end{equation}
Since we have
\begin{equation*}
\begin{split}
&\frac{d g_1^{n-1}(u)}{d u} = \\
&\frac{1}{2 \alpha_1}\left(-2\alpha_1 + \frac{4\alpha_1 (\beta_1^{n-1} - \alpha_1 + A^{n-1})}{2 \sqrt{(\beta_1^{n-1} - 2\alpha_1 u)^2 - 4 \alpha_1 (\gamma_1^{n-1} + u (A^{n-1} + \alpha_1) + \alpha_1 u^2)}}\right) = -1,
\end{split}
\end{equation*}
then (\ref{der_coeff_w0}) under the condition (\ref{CFL}) reduces to $k \le 1.$

 
Reasoning as above, for arc $m$ we obtain:
\begin{equation}\label{u26+_boundary:2}
\begin{split}
u_{+,m}^{n+1,M_{m}+1}&= \left(1-\frac{\lambda_{m}k}{h_{m}}-\frac{k}{4} -\frac{k}{4\lambda_{m}}\chi \phi_{x,m}^{n,M_{m}+1}\right) u_{+,m}^{n,M_{m}+1} \\
& + g_2^{n-1}(u_{+,m}^{n,M_{m}+1})\left(\frac{k}{4}+\frac{k}{4\lambda_{m}}\chi \phi_{x,m}^{n,M_{m}+1}\right)\\
&+  u_{-,m}^{n,M_{m}}\left(\frac{k}{4} + \frac{k \chi}{4 \lambda_{m}}\frac{\phi_{m}^{n,M_{m}+1}-\phi_{m}^{n,M_{m}-1}}{2 h_{m}} \right)\\
&+ u_{+,m}^{n,M_{m}}\left(\frac{\lambda_{m} k}{h_{m}}-\frac{k}{4} + \frac{k \chi}{4 \lambda_{m}}\frac{\phi_{m}^{n,M_{m}+1}-\phi_{m}^{n,M_{m}-1}}{2 h_{m}}\right).
\end{split}
\end{equation}

As above, the condition (\ref{CFL}) reduces to $k \le 1.$



\subsection{Approximation of the parabolic equation for $\phi$}

We  solve the parabolic equation, using a finite differences scheme in space and a Crank-Nicolson method in time, namely an explicit-implicit method in time. 

Therefore, we will have the following equation for $\phi^{n+1,j}_{i}, \, 1 \leq j \leq M_{i}$,
\begin{equation}\label{discr-CN}
\begin{aligned}
\phi^{n+1,j}_{i}&=\phi^{n,j}_{i}-\frac{D_{i}k}{2h_{i}^2} \left( -\phi^{n,j+1}_{i}+2\phi^{n,j}_{i}-\phi^{n,j-1}_{i}\right)\\
&-\frac{D_{i}k}{2h_{i}^2} \Bigl( -\phi^{n+1,j+1}_{i}+2\phi^{n+1,j}_{i}-\phi^{n+1,j-1}_{i}\Bigr)\\
&+\frac{a_{i}k}{2}(u^{n+1,j}_{i}+u^{n,j}_{i})-\frac{b_{i}k}{2}(\phi^{n+1,j}_{i}+\phi^{n,j}_{i}).
\end{aligned}
\end{equation}
Now, let us find the two boundary conditions needed on each interval. As in \cite{BNR14}, the boundary conditions will be given in the case of an outer node  and in the case of an inner node.
On the outer boundary we discretize $\phi_x$ using a second order approximation. Then, the numerical condition, corresponding to (\ref{general_bound_phi}), for a general $\bar\phi$, not necessarily depending on $\phi_{i}^{n+1}$, reads as:
  \begin{equation}\label{discr-outer_nonull}
\left\{\begin{aligned}
&\phi_{i}^{n+1,0}= \frac 4 3 \phi^{n+1,1}_{i}-\frac 1 3\phi^{n+1,2}_{i} - \frac{2 h_i}{3}\bar\phi, \, \textrm{ if } i\in I_{out}, \ \\
&\phi_{i}^{n+1,M_{i}+1}= \frac 4 3 \phi^{n+1,M_{i}}_{i}-\frac 1 3\phi^{n+1,M_{i}-1}_{i}+ \frac{2  h_i}{3}\bar\phi, \, \textrm{ if } i\in O_{out}.
\end{aligned}\right.
\end{equation}

Let us now describe our numerical approximation for the transmission condition  \eqref{transmission-phi} which,  as the transmission condition  for the hyperbolic part  \eqref{transmission}, couples the $\phi$ functions of arcs having  a node in common.
 
At  node $p$ condition \eqref{transmission-phi} is discretized using the second-order discretization as described in \cite{BNR14}:
\begin{equation}\label{discr-transmi}
\begin{aligned}
 \eta_{i}^p \phi_{i}^{n+1,M_{i}+1}=&\frac{4}{3}\phi_{i}^{n+1,M_{i}}-\frac{1}{3}\phi_{i}^{n+1,M_{i}-1} 
+\frac{2}{3}\frac{h_{i}}{D_{i}}\sum_{j \in I_{p}} \kappa_{i,j} \phi_{j}^{n+1,M_{j}+1}
\\& \quad \quad  \qquad  \qquad 
+\frac{2}{3}\frac{h_{i}}{D_{i}}\sum_{j \in O_{p}} \kappa_{i,j} \phi_{j}^{n+1,0} ,  
  \textrm{ if }i\in I_{p}, 
\\
 \eta_{i}^p \phi_{i}^{n+1,0}=&\frac{4}{3}\phi_{i}^{n+1,1}-\frac{1}{3}\phi_{i}^{n+1,2}
+\frac{2}{3}\frac{h_{i}}{D_{i}}\sum_{j \in I_{p}} \kappa_{i,j} \phi_{j}^{n+1,M_{j}+1}
\\& \quad \quad  \qquad  \qquad 
+\frac{2}{3}\frac{h_{i}}{D_{i}}\sum_{j \in O_{p}} \kappa_{i,j} \phi_{j}^{n+1,0},   
\textrm{ if }i\in O_{p},
\end{aligned}
\end{equation}
with $ \eta_{i}^p= 1+\frac{2}{3}\frac{h_{i}}{D_{i}} \sum_{j\in I_{p} {\cup} O_{p}} \kappa_{i,j}$.
Let us remark that the previous discretizations are compatible with relations \eqref{discr-outer_nonull} with $\bar\phi=0$ (homogeneous Neumann boundary condtions) considering that for outer boundaries the coefficients $\kappa_{i,j}$ are null. Therefore, in this case, the value of $\eta_{i}^{out}$  is just equal to  $1$.
Since equations  \eqref{discr-transmi} are coupling of the unknowns of all arcs altogether, we have to solve  a large system which contains  all the equations of type \eqref{discr-CN} and also  the discretizations of transmission conditions \eqref{discr-transmi}.
Note that for the computational resolution of the mentioned system, characterized by a sparse banded  matrix, we used the LAPACK-Linear Algebra PACKage routine DGBSV designed for banded matrix.
Once the values of $\ds \phi^{n+1,j}_{i}$ are known, we can compute a  second-order discretization of the derivatives of $\phi$  which gives the values of the $f$ function, namely :
\begin{equation*}
\phi_{x,i}^{n+1,j}=  \left\{\begin{aligned}
& \ds\frac{1}{2\,h_{i}} \left(\phi ^{n+1,j+1}_{i}-\phi ^{n+1,j-1}_{i}\right),  1 \leq j \leq M_{i},\\
& \ds\frac{1}{2\,h_{i}} \left( -\phi ^{n+1,2}_{i}+4\phi ^{n+1,1}_{i}-3\phi ^{n+1,0}_{i}\right), j=0,\\
& \ds\frac{1}{2\,h_{i}} \left(\phi ^{n+1,M_{i}-1}_{i}-4\phi ^{n+1,M_{i}}_{i}+3\phi ^{n+1,M_{i}+1}_{i}\right), j=M_{i}+1.
\end{aligned}\right.
\end{equation*}
The discretization of $f$ needed at equations \eqref{scheme-uv-interval},\eqref{bound3}, and \eqref{bound4} is therefore given by $f^{n+1,j}_{i}= \chi \phi_{x,i}^{n+1,j}  u^{n+1,j}_{i}$.

\section{Tests on different network geometries} \label{sec4}
Here we present some computational tests on different network topologies to investigate the behaviour of individuals moving through them. We firstly consider the basic network composed of one incoming and two outgoing arcs (T-shaped network) in order to show that the model is able to reproduce the main features of the plasmodium. Then we consider a network of 7 arcs and 6 nodes (diamond-like graph) and a more complex network of 26 arcs and 18 nodes, both of them connecting two exits (a source and a sink), where inflow conditions are implemented to mimic the presence of food sources. Finally, we consider a network of 21 arcs and 15 nodes with multiple exits (two sources and three sinks).\\
For all the networks we deal with we assume to have equal velocities $\lambda_i=\lambda$ and equally distributed transmission coefficients at each node. If, for instance, we have a node $p$ connected with the total number of $N$ arcs (no matter how many are incoming or outgoing), we assume to have:
\begin{equation}\label{cond_xi}
\xi_{ij}=\frac{1}{N}, \forall i \in I_p, \forall j \in O_p, 
\end{equation}
in such a way to satisfy the conditions \eqref{condition_lambda} and \eqref{dissip-coeff}.

\subsection{T-shaped graph: four nodes and three edges}
\begin{figure}[htbp!]
\begin{center}
\includegraphics[width=5.5cm]{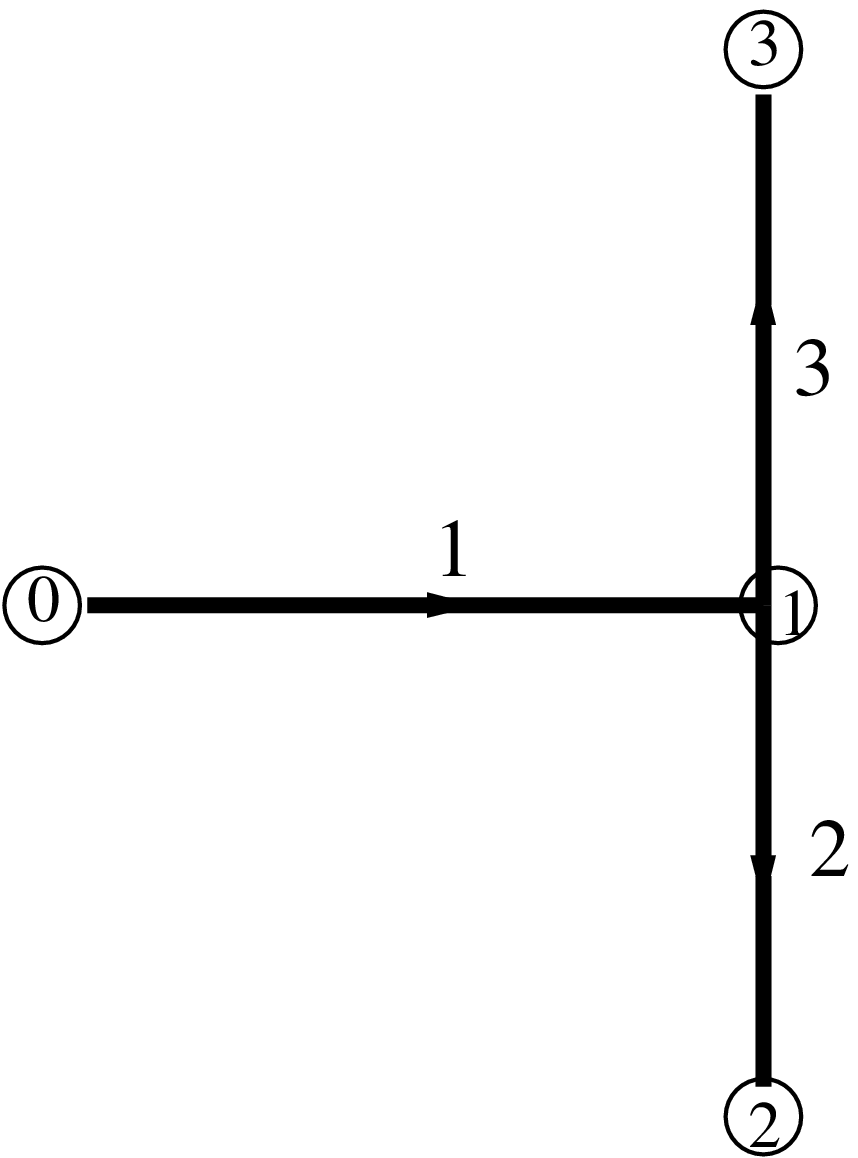}
\caption{Network of three edges: one incoming and two outgoing arcs (T-shaped network).}
\label{fig:1}
\end{center}
\end{figure}

Here we consider the network of three arcs and four nodes (1 internal node and three external nodes) shown in Figure \ref{fig:1}. We assume to have arcs of length $L_i=1$ and we set parameters as $a_i=1, b_i=0.1, \lambda_i=\sqrt{0.33}, D_i=1$, with $\chi_i= 1$ representing positive chemotaxis, for $i=1,2,3$. Furthermore, for the transmission conditions for $u$ at the internal node we set dissipative coefficients $\xi_{i,j}=\frac{1}{3}$ for $i,j=1,2,3$ and for the transmission conditions for $\phi$ we assume $\kappa_{i,j}=1$ for $i \neq j$ and $\kappa_{i,i}=0$, for $i,j=1,2,3$.

\subsubsection{Example 1: the zero flux case}
\begin{figure}[htbp!]
\begin{center}
\includegraphics[height=5.5cm,width=5.5cm]{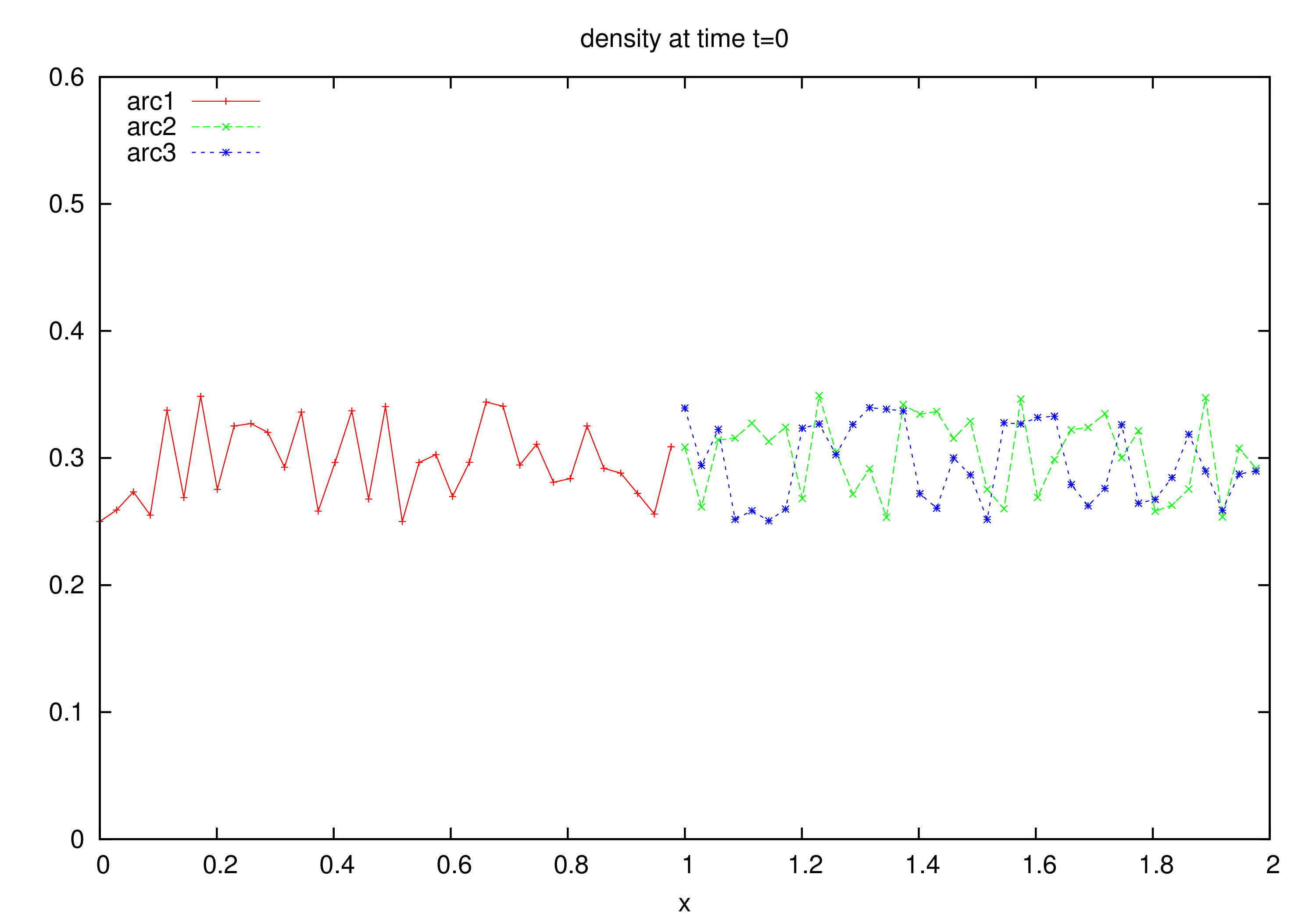}\quad \includegraphics[height=5.5cm,width=5.5cm]{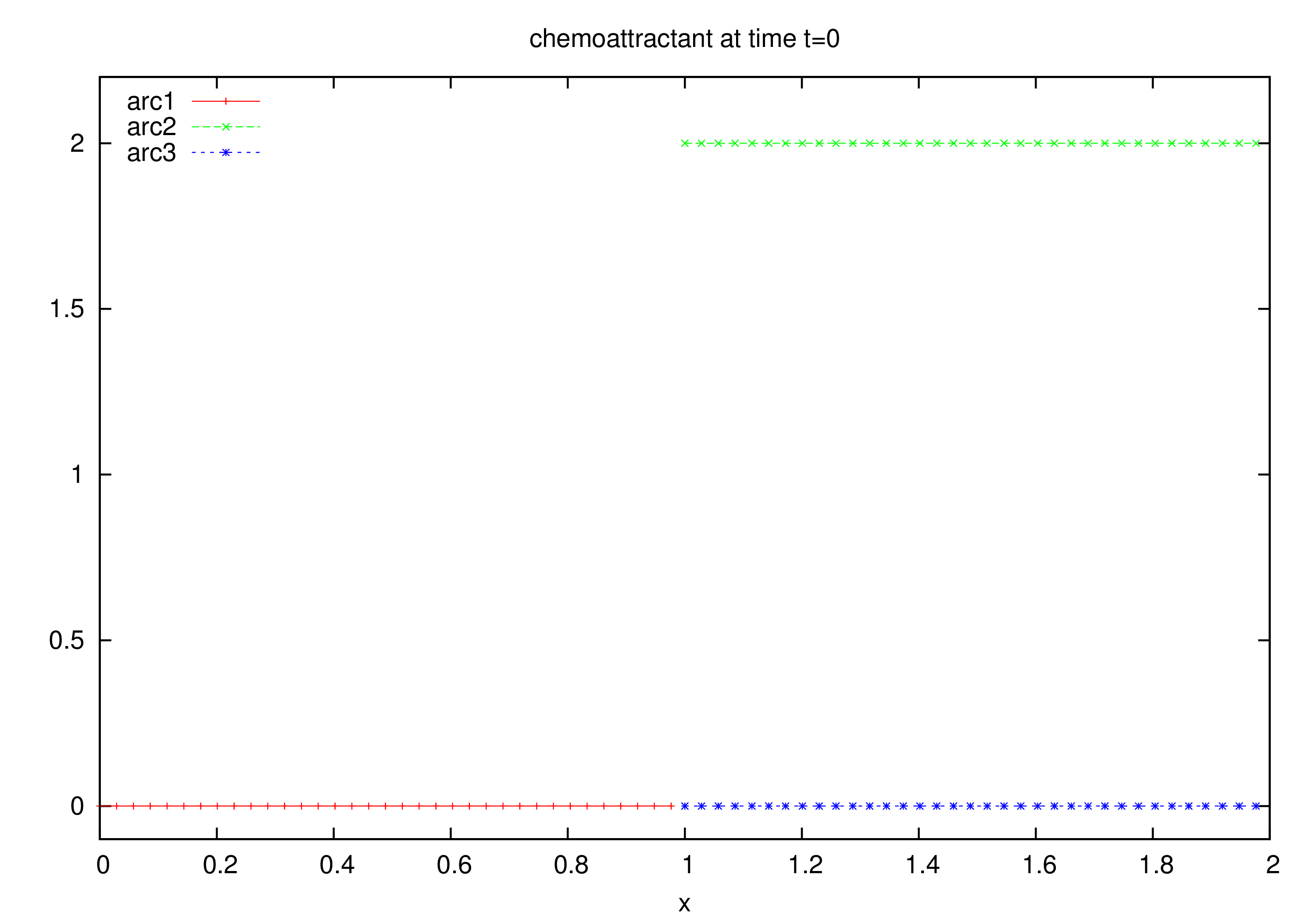}
\caption{Initial data for the test in Example 1.}
\label{fig:data}
\end{center}
\end{figure}

In this first example we set homogeneous Neumann boundary conditions at the outer boundaries for both $u$ and $\phi$, in order to have zero flux conditions as in \cite{BGKS}. The initial values for $u_i$ are randomly equally distributed in $[0.25,0.35]$ for each $i=1,2,3$. Moreover, we set $\phi_1 (x, 0) = \phi_3 (x, 0) = 0$ and $\phi_2 (x, 0) = 2$, see Fig. \ref{fig:data}.

\begin{figure}[htbp!]
\begin{center}
\includegraphics[height=5.5cm,width=5.5cm]{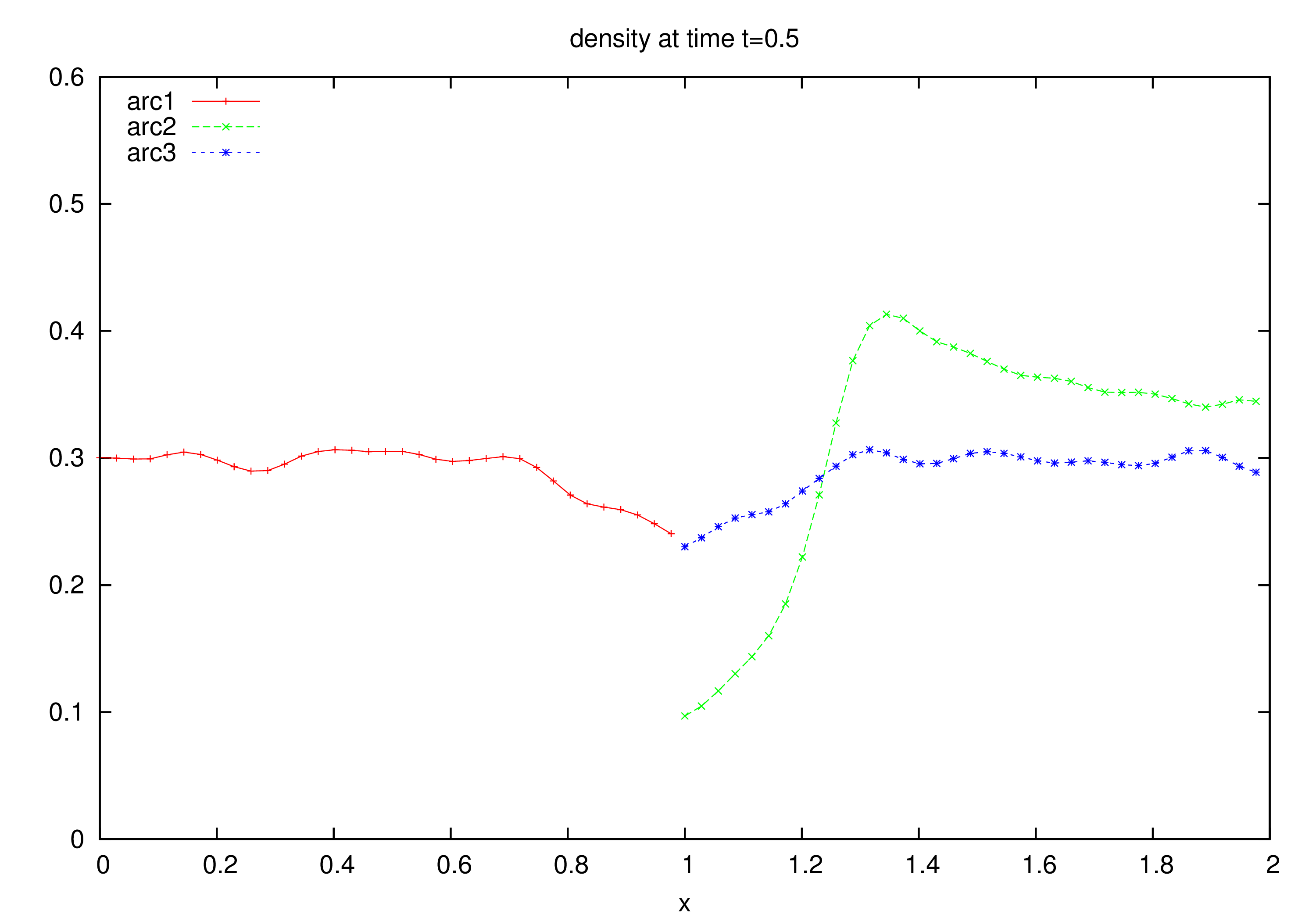}\quad \includegraphics[height=5.5cm,width=5.5cm]{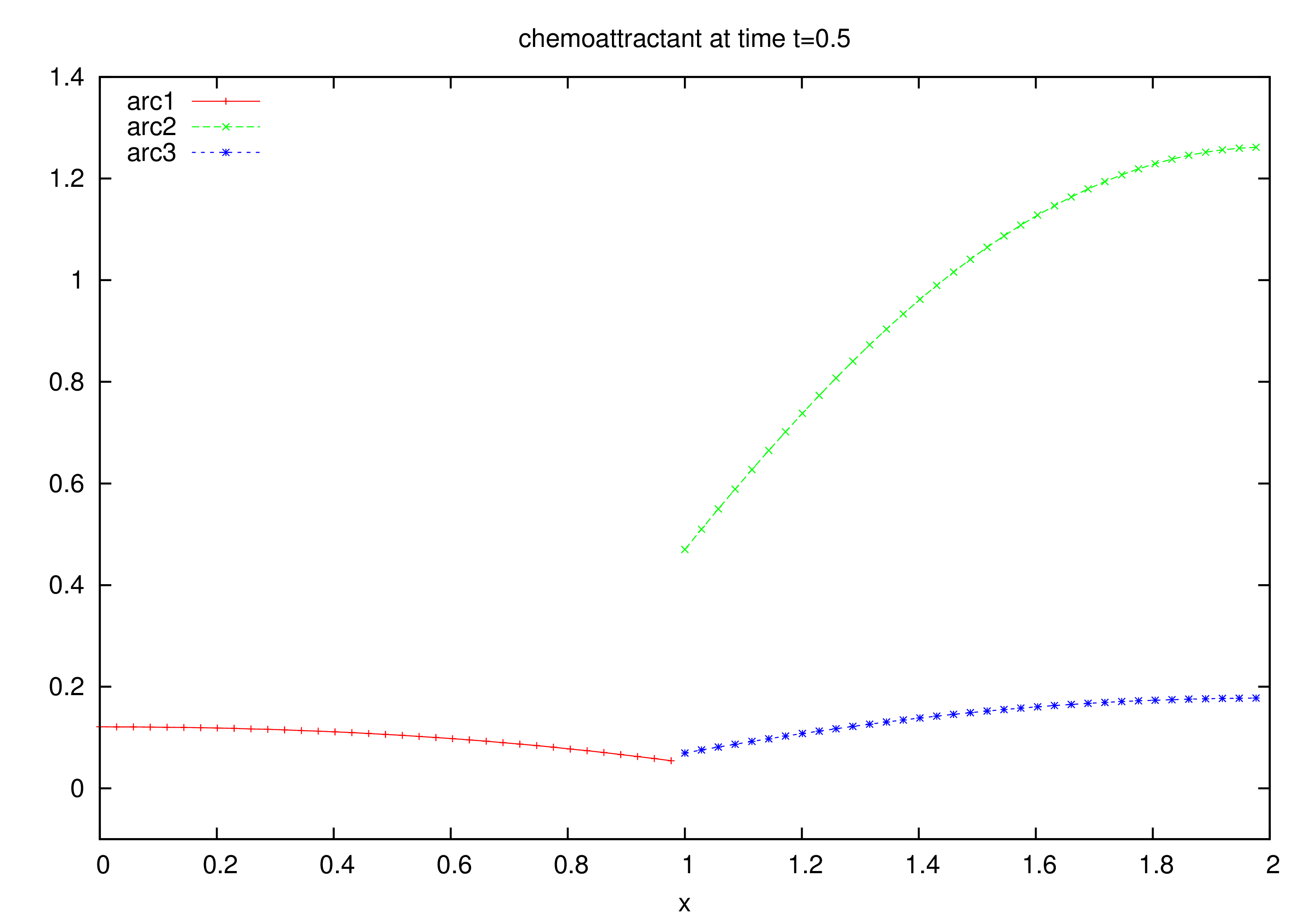}
\caption{Example 1. The distribution of the density $u_i(x,t)$ (on the left) and of the chemoattractant $\phi_i(x,t)$ (on the right), for $i = 1, 2, 3$, inside the T-shaped graph at time $t=0.5$.}
\label{fig:2}
\end{center}
\end{figure}

\begin{figure}[htbp!]
\begin{center}
\includegraphics[height=5.5cm,width=5.5cm]{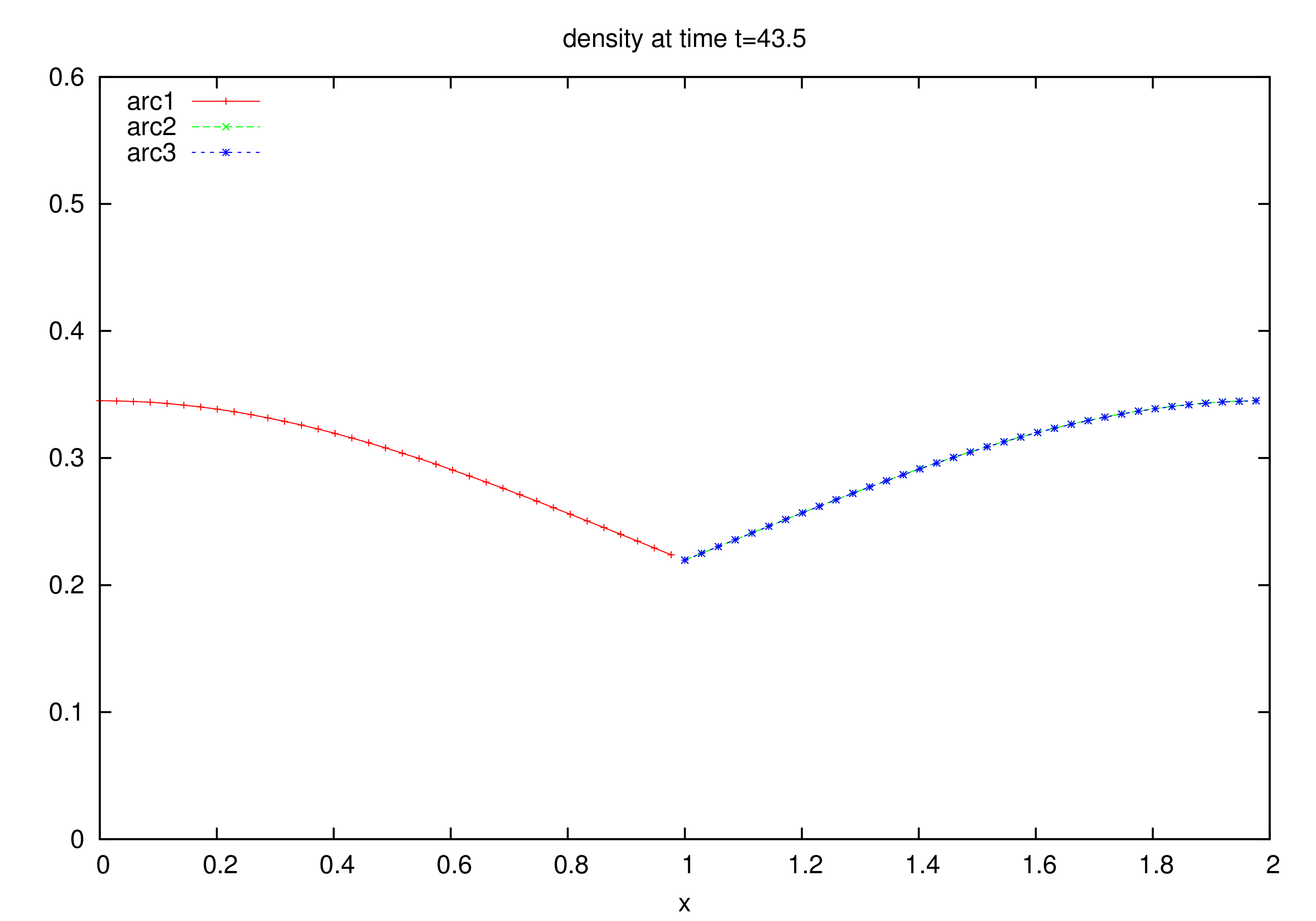}\quad \includegraphics[height=5.5cm,width=5.5cm]{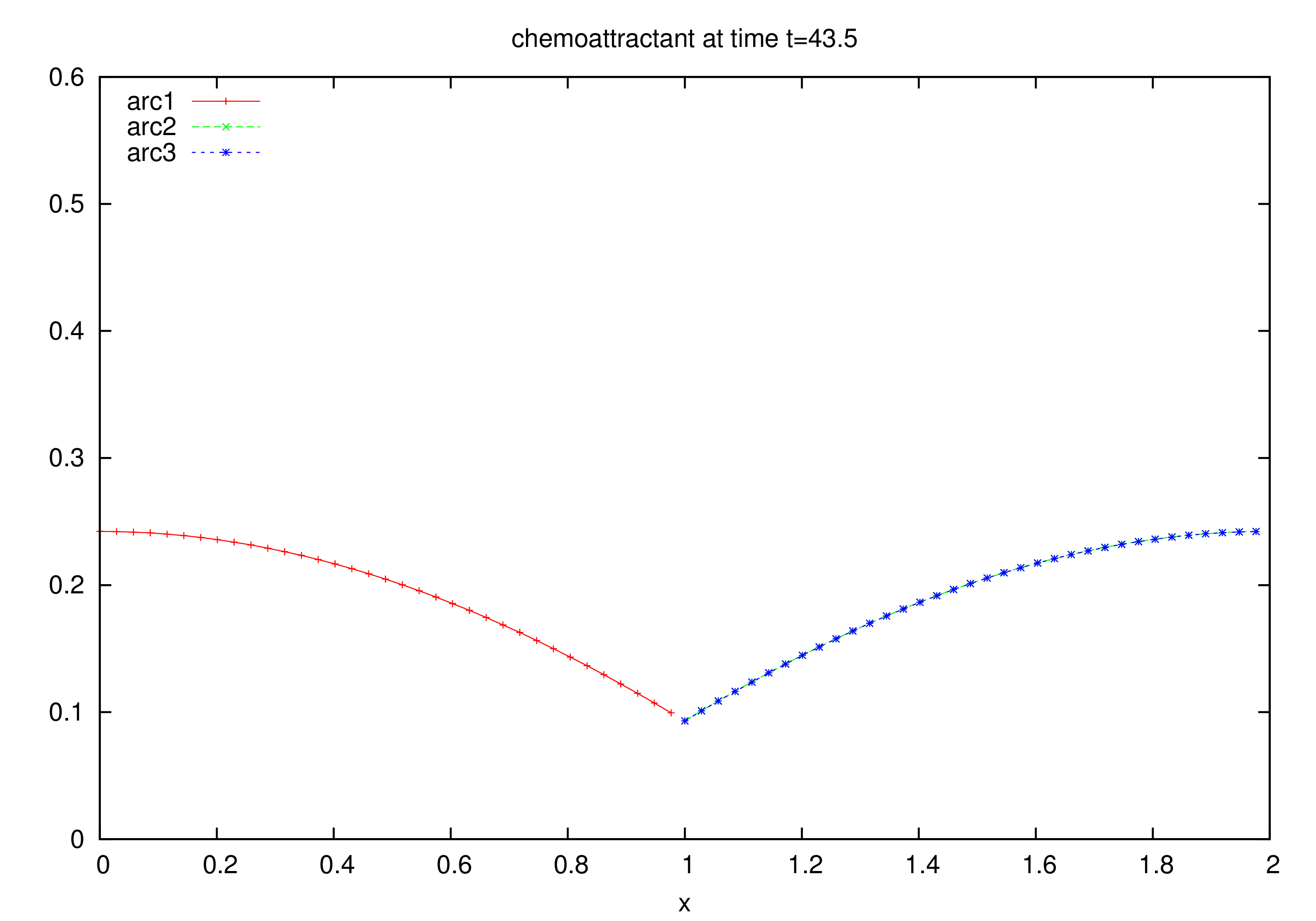}
\caption{Example 1. The distribution of the density $u_i(x,t)$, (on the left) and of the chemoattractant $\phi_i(x,t)$ (on the right), $i = 1, 2, 3$, inside the  T-shaped graph at time $t=43.5$.}
\label{fig:4}
\end{center}
\end{figure}

\begin{figure}[htbp!]
\begin{center}
\includegraphics[height=7cm,width=7cm]{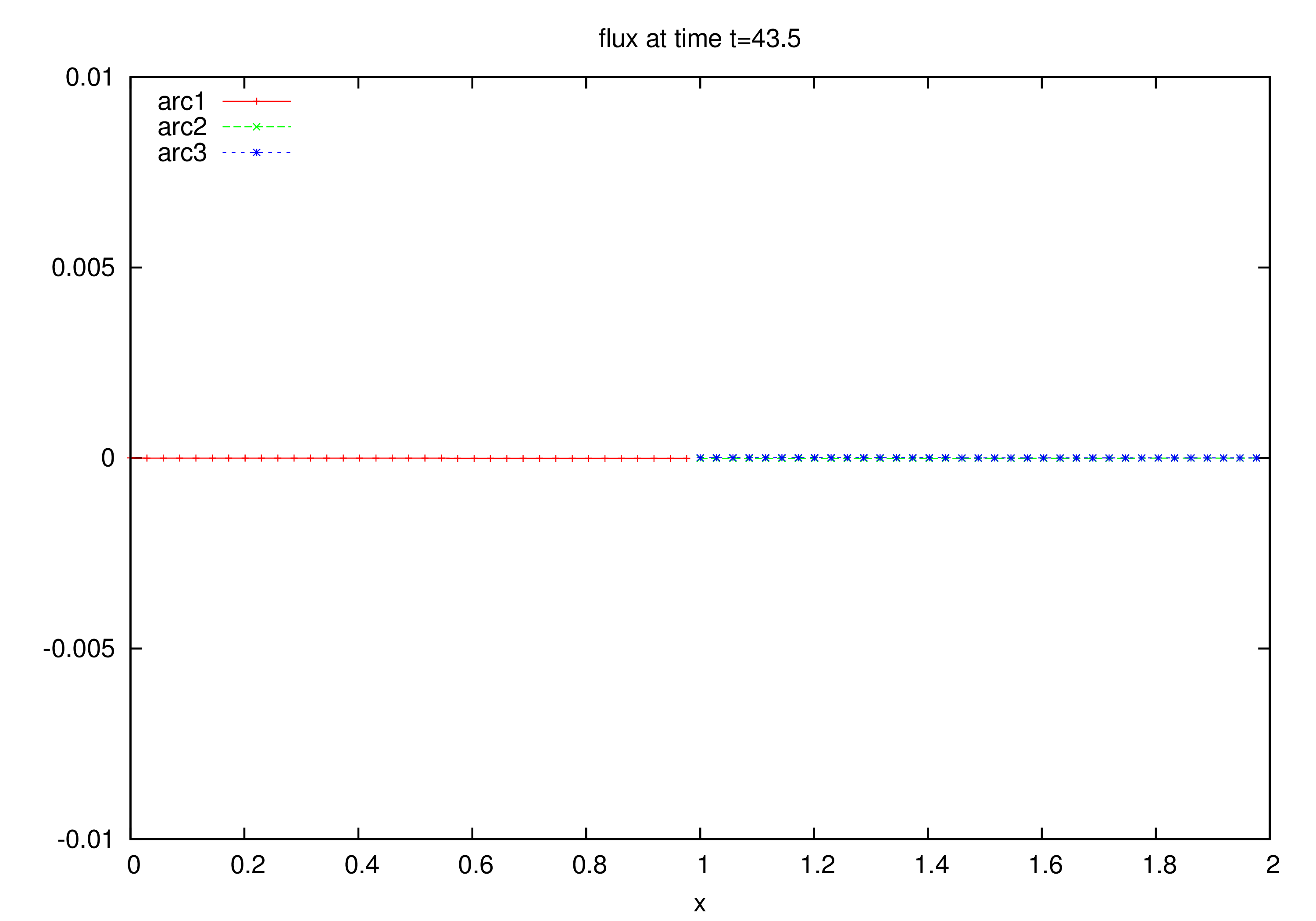}
\caption{Example 1. The asymptotic fluxes of the arcs of the T-shaped graph at time $t=43.5$.}
\label{fig:5}
\end{center}
\end{figure}

As shown in Fig. \ref{fig:1}, the non-zero initial condition for chemoattractant on arc $2$, leads the cells to move towards it. Since there the concentration of chemoattractant is largest, they accumulate at the right boundary of the arc, see the left picture of Fig. \ref{fig:2}. 
Moreover, due to the diffusion flux, chemoattractant distributes on the other arcs, thus determining a fast decrease in the quantity of the chemoattractant on arc $2$ as shown in the right picture of Fig. \ref{fig:2}. The mentioned results are in accordance to the ones reported in \cite{BGKS}.
We also computed the asymptotic solution achieved at time $T=43.5$ and in Figg. \ref{fig:4}-\ref{fig:5} we reported the profiles of the final densities and fluxes.

\begin{figure}[htbp!]
\begin{center}
\includegraphics[height=5.5cm,width=5.5cm]{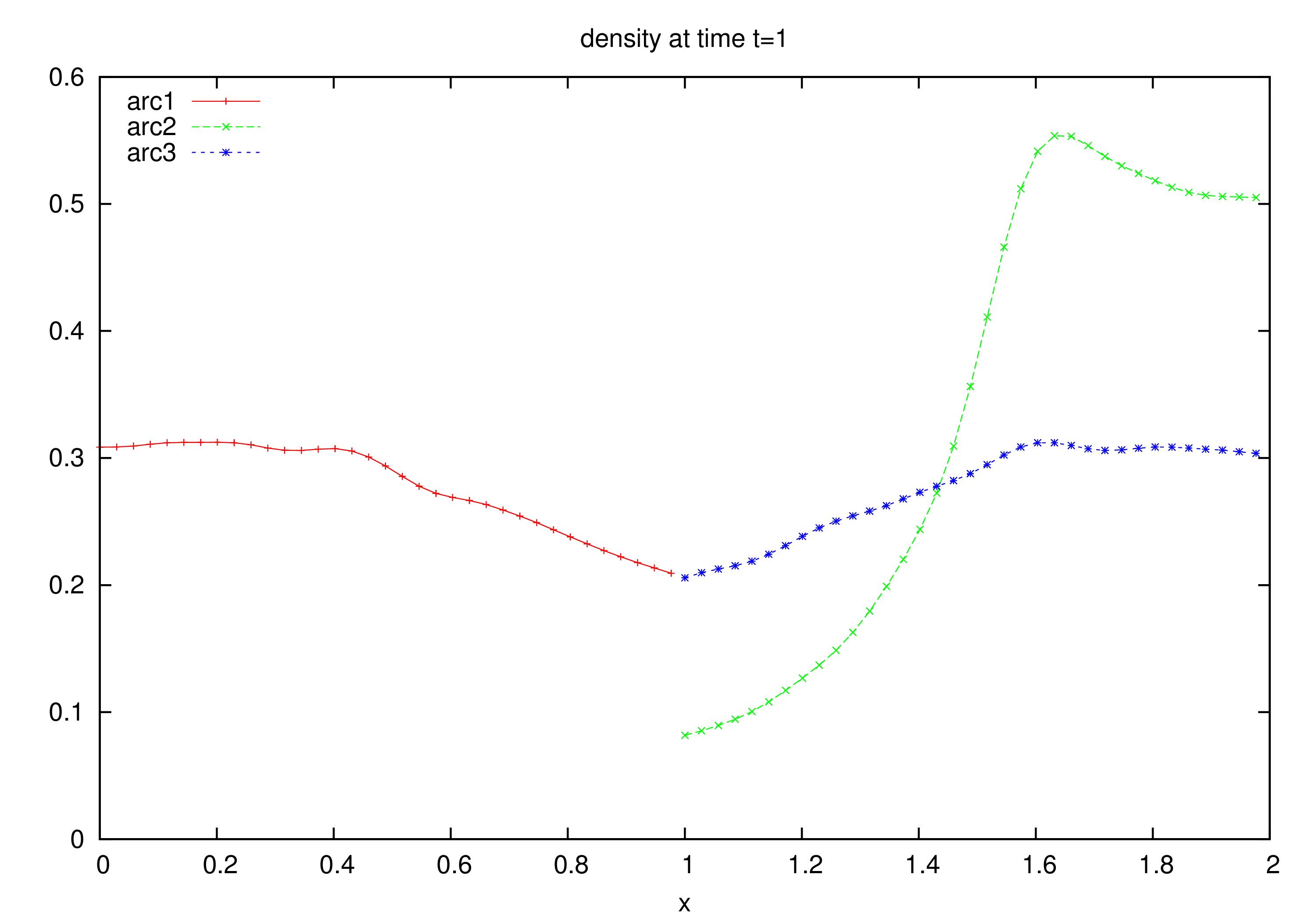}\quad \includegraphics[height=5.5cm,width=5.5cm]{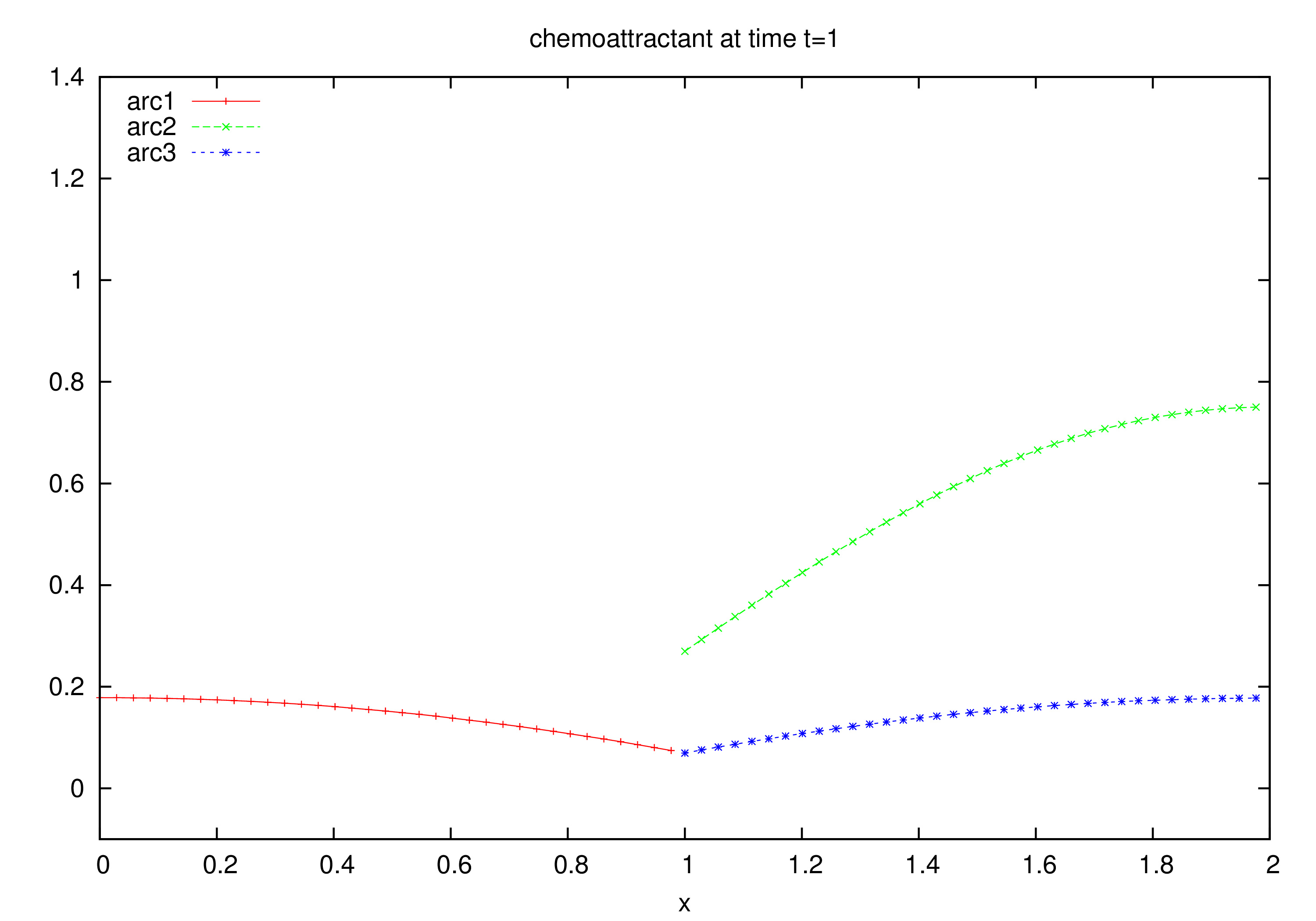}
\caption{Example 1. The distribution of the density $u_i(x,t)$, (on the left) and of the chemoattractant $\phi_i(x,t)$ (on the right), $i = 1, 2, 3$, inside the T-shaped graph at time $t=1$.}
\label{fig:3}
\end{center}
\end{figure}

\subsubsection{Example 2: non-homogeneous Neumann boundary conditions}

Now, in order to reproduce the experiment with the T-shaped plasmodium tubes of Physarum presented in \cite{KTN2007} to show the property of dead-end cutting, we set non-homogeneous Neumann boundary conditions at the outer boundaries for $\phi$, to mimic the presence of food at the left and lower ends.  For $u$ we impose homogeneous Neumann boundary conditions at the outer boundaries as above. The initial values for $u_i$ are randomly equally distributed in $[0.25,0.35]$ for each $i=1,2,3$ to mimic the plasmodium spread in the network, while we set $\phi_i (x, 0) = 0, i=1,2,3$.

\begin{figure}[htbp!]
\begin{center}
\includegraphics[height=5.5cm,width=5.5cm]{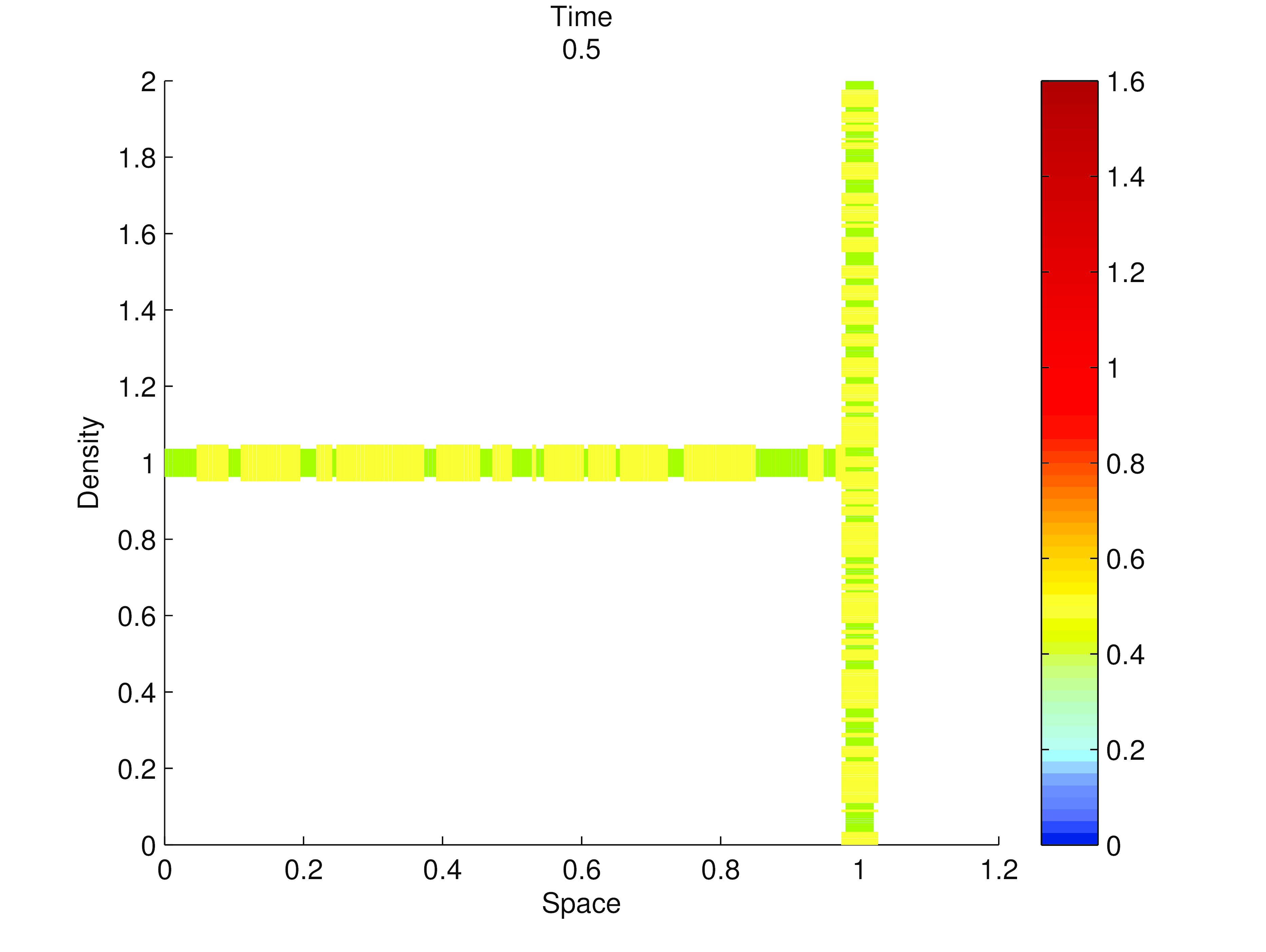}\quad \includegraphics[height=5.5cm,width=5.5cm]{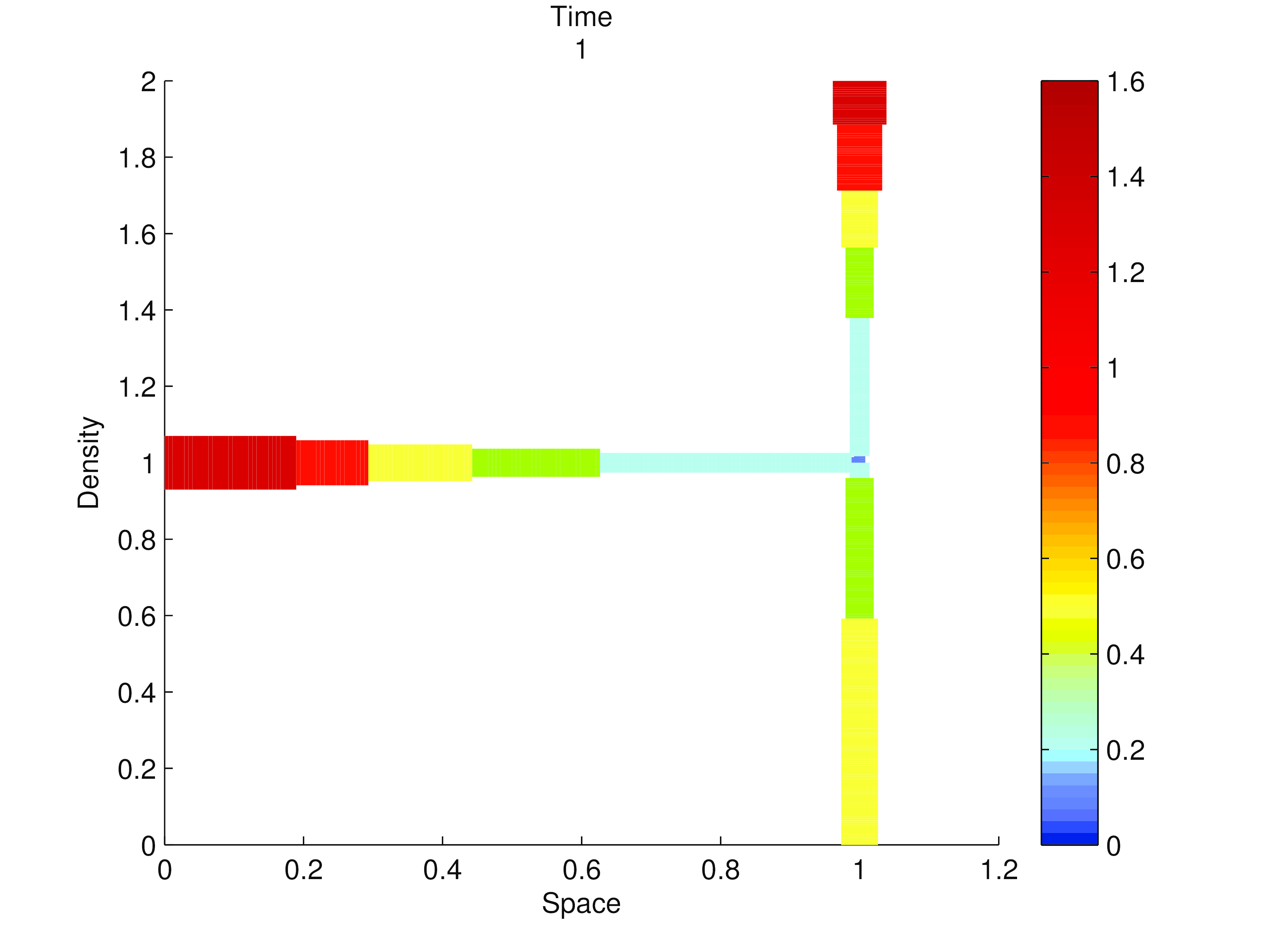}
\caption{Example 2. The distribution of the density $u_i(x,t)$ on each arc of the T-shaped network at time $t=0.5$ (on the left) and $t=1$ (on the right).}
\label{fig:1t}
\end{center}
\end{figure}

\begin{figure}[htbp!]
\begin{center}
\includegraphics[height=5.5cm,width=5.5cm]{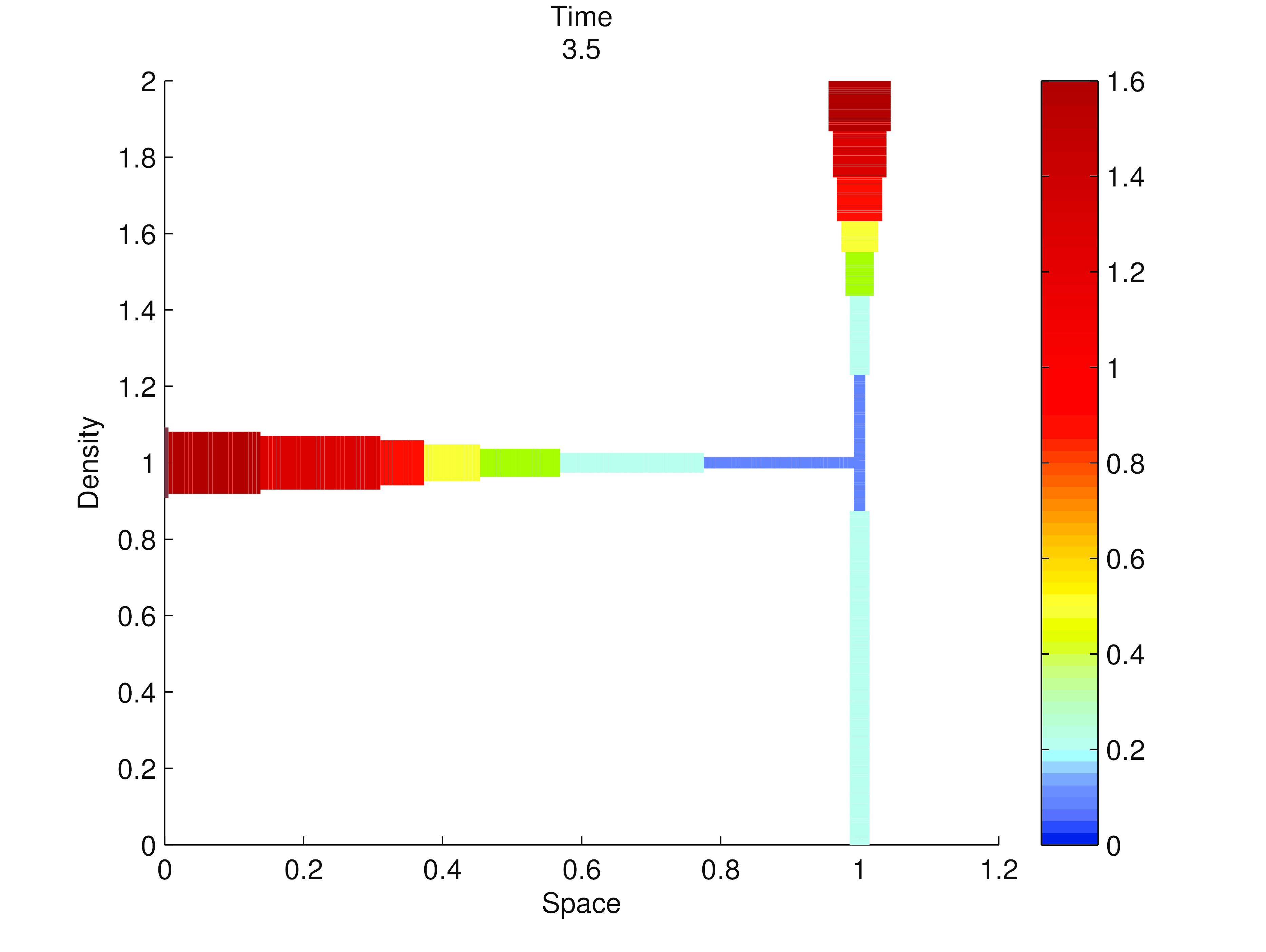}\quad \includegraphics[height=5.5cm,width=5.5cm]{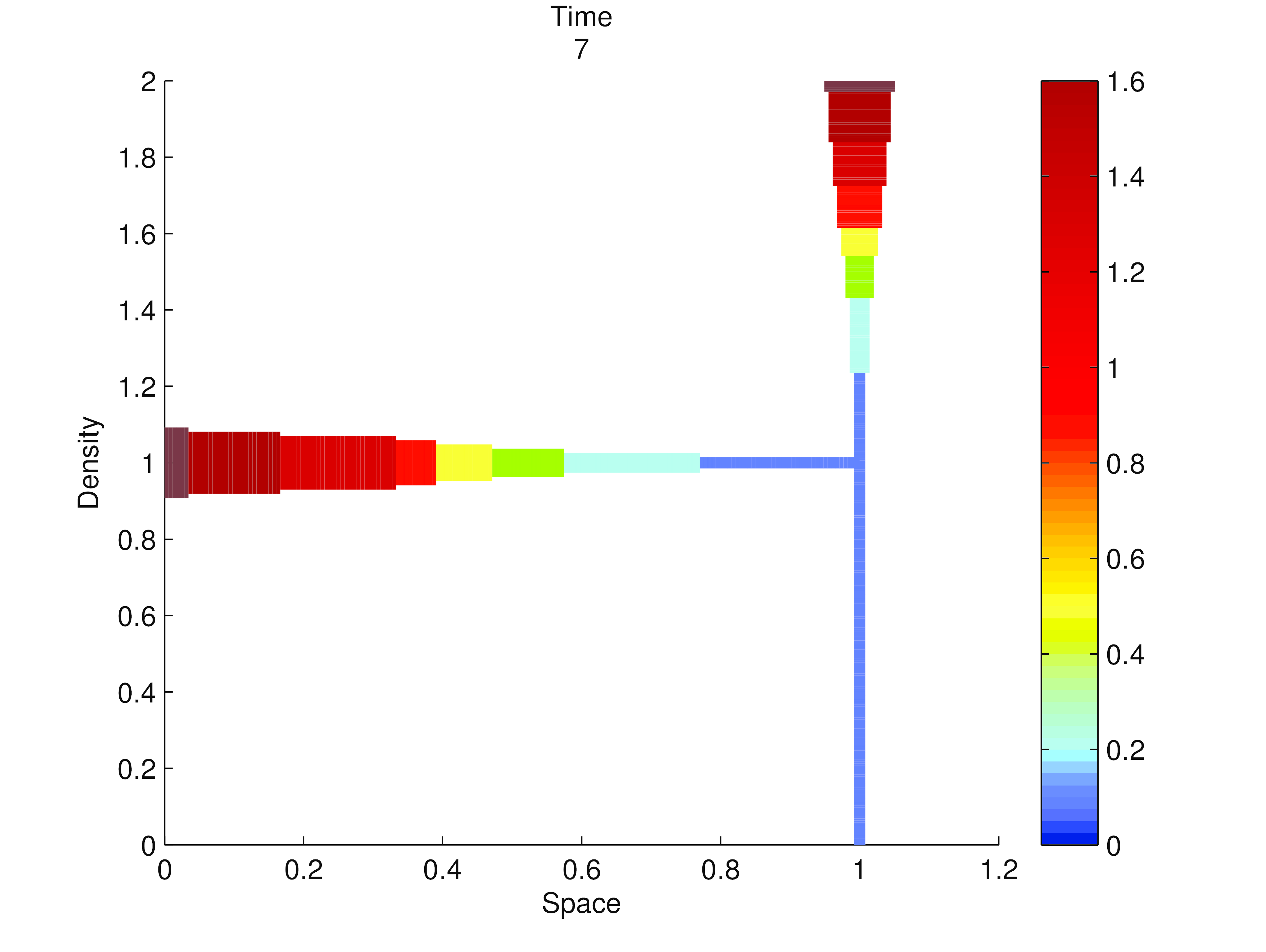}
\caption{Example 2. The distribution of the density $u_i(x,t)$ on each arc of the T-shaped network at time $t=3.5$ (on the left) and $t=7$ (on the right).}
\label{fig:2t}
\end{center}
\end{figure}

Let us recall that in our modeling the dead-end cutting property corresponds to have a mass equal to zero on the egdes where the flux is null. In accordance to the cited experiment, as the organism starts moving to feed itself, we observe that the mass concentrate on the left and lower edges (arcs 1 and 2) since food sources are located at the outer nodes $0$ and $2$, and decreases quickly until it becomes null on the upper edge (arc 3), whose endpoint has no food source.

\subsection{Study on the path followed by amoeboid-like organism.}
Here we are interested in the identification of the shortest path followed by individuals in the case of positive (attractive) chemotaxis. 
\subsubsection{Wheatstone bridge-shaped network with a source and a sink (two exits)}
\begin{figure}[htbp!]
\begin{center}
\includegraphics[width=8cm]{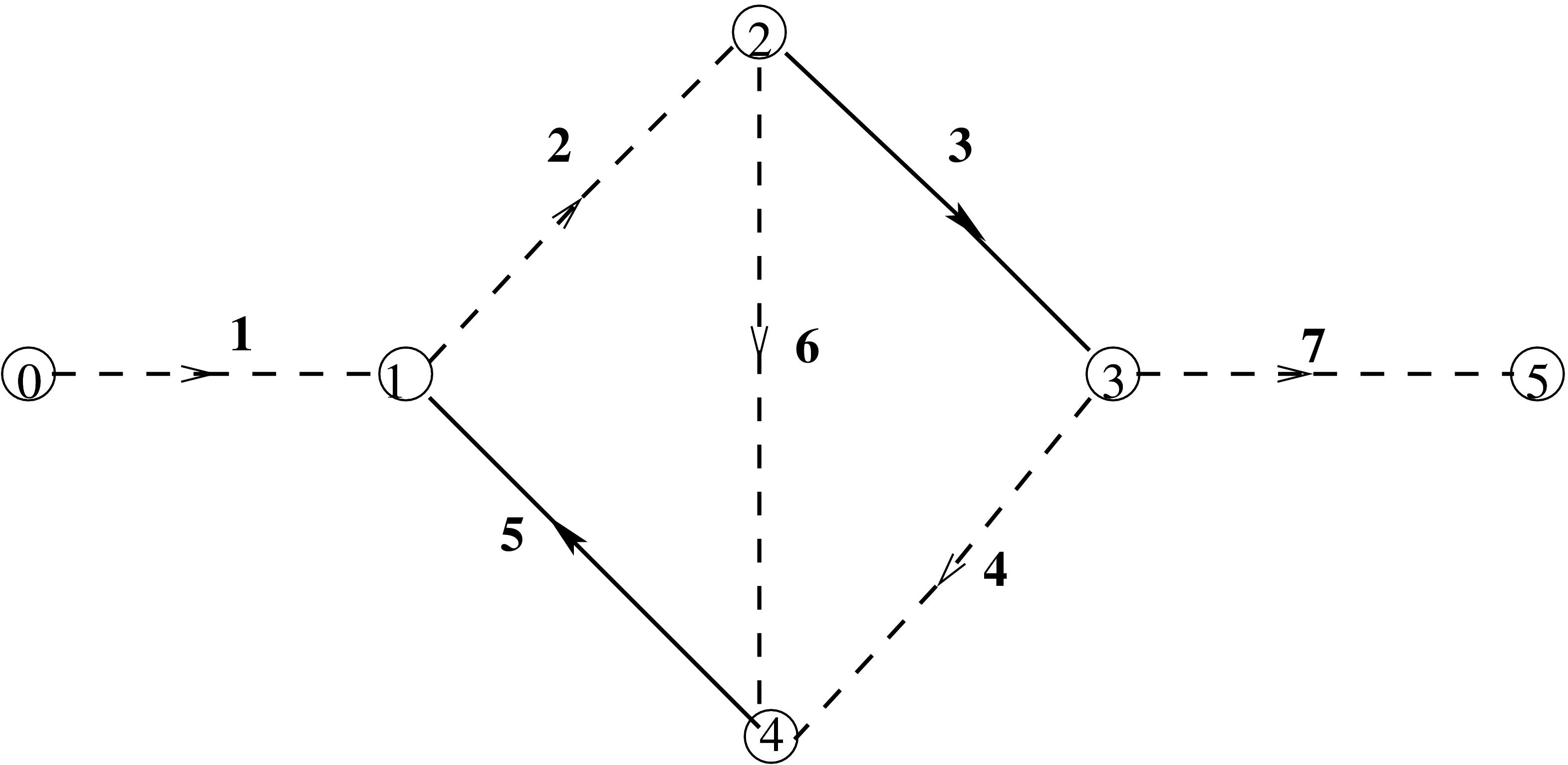}
\caption{Wheatstone bridge-shaped network: a network with 7 arcs and 4 internal nodes.}
\label{diamondgraph}
\end{center}
\end{figure}

Let us consider the Wheatstone bridge-shaped network where, for our convenience, we inserted one incoming and one outgoing arc, respectively, exiting or entering an external node, as shown in Fig. \ref{diamondgraph}. Since this modification does not change the structure of the shortest paths in the network, the analysis on the Wheatstone bridge network can be found in \cite{Mi2007}, where it was proven that the globally asymptotically stable equilibrium point corresponds to the shortest path
connecting the exits.
If, for instance, we want the path $1 \to 2 \to 6 \to 4 \to 7 $  to be shortest one, we need to set:
\begin{equation}
L_{2} + L_{6} < L_{5}  \ \textrm{and} \  L_{4} + L_{6} < L_{3}.
\end{equation}
In particular, we assume 
\begin{displaymath}
L_{1} = L_7 = 0.2,\ L_{2} = L_4 = 0.3, \ L_{3} = L_{5} = 2,  L_{6} = 0.3.
\end{displaymath}
Note that in Fig. \ref{diamondgraph} we depicted the arcs composing the shortest path with the dashed line.\\
We set parameters as $a_i=1, b_i=0.1, \lambda_i=\sqrt{0.33}, \chi_i= 1$, $D_i=1$, for $i=1,\ldots,7$. As before, for the transmission conditions for $u$, at each node $p$ we set dissipative coefficients $\xi_{i,j}=\frac{1}{3}$ for $i,j=1,\ldots,7$ and for the transmission conditions for $\phi$ we assume $\kappa_{i,j}=1$ for $i \neq j$ and $\kappa_{i,i}=0$, for $i,j=1,\ldots,7$.
 The initial values for $u_i$ are randomly equally distributed in $[0.45,0.55]$ and we set $\phi_i (x, 0) =0$ for each $i=1,\ldots,7$.

We set non-homogeneous Neumann boundary conditions for the chemottractant $\phi$ at the node $0$ and at the node $5$ (inflow conditions) :
\begin{equation}\label{inflow_cond_diamond} 
\partial_x \phi_1(0,t) = -1, \quad \partial_x \phi_{7}(L_7,t) = 1.
\end{equation}

Firstly, we assume for $u$ homogeneous Neumann boundary conditions $u_x(\cdot,t)=0$ corresponding to $v(\cdot,t)=0$ at the node $0$ and at the node $5$.

\begin{figure}[htbp!]
\begin{center}
\includegraphics[height=5.5cm,width=5.5cm]{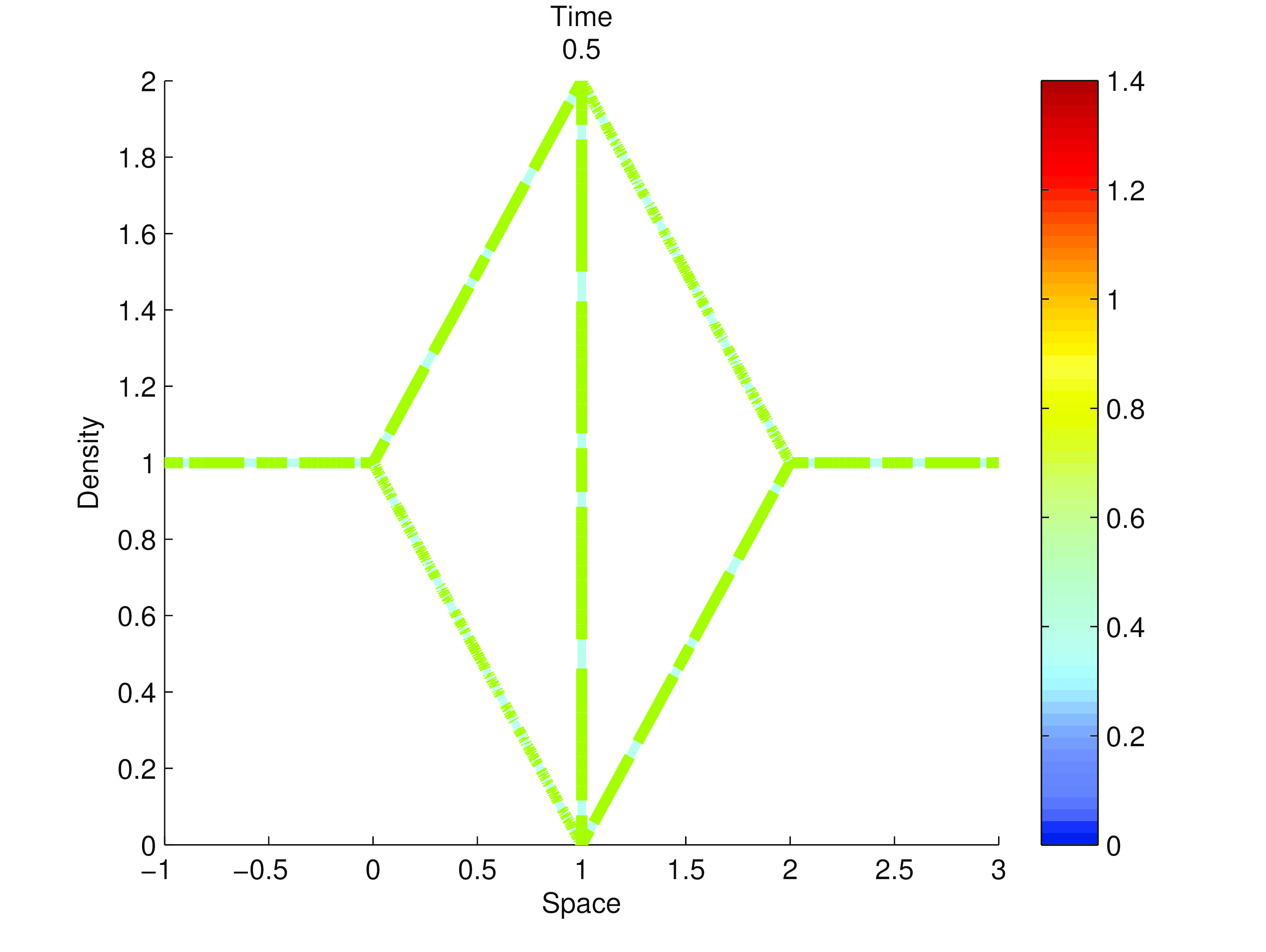}\quad \includegraphics[height=5.5cm,width=5.5cm]{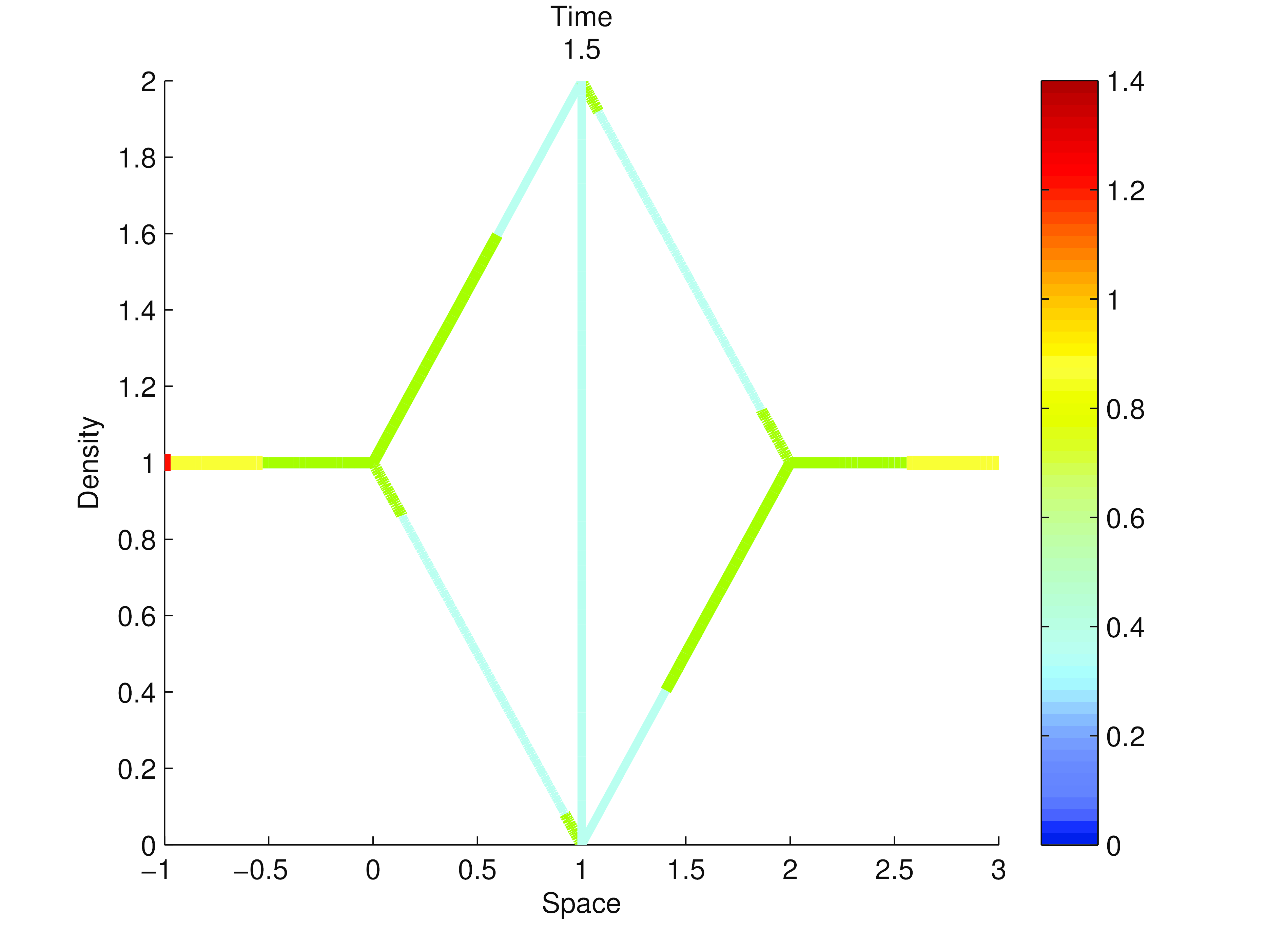}
\caption{Wheatstone bridge-shaped network. The distribution of the density $u_i(x,t)$ on each arc of the diamond graph-like network at time $0.5$ (on the left) and $t=1.5$ (on the right).}
\label{fig:1d}
\end{center}
\end{figure}

\begin{figure}[htbp!]
\begin{center}
\includegraphics[height=5.5cm,width=5.5cm]{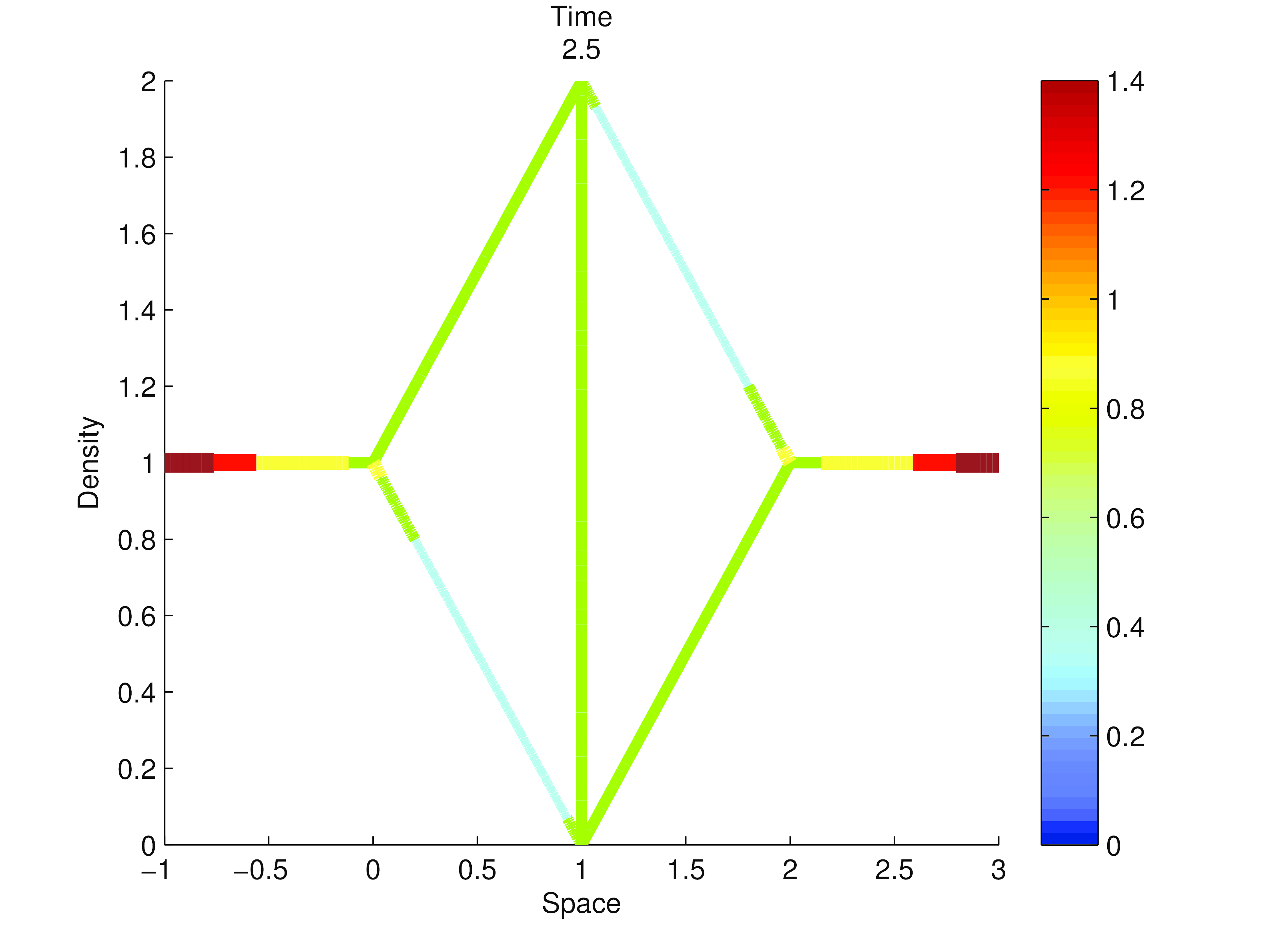}\quad \includegraphics[height=5.5cm,width=5.5cm]{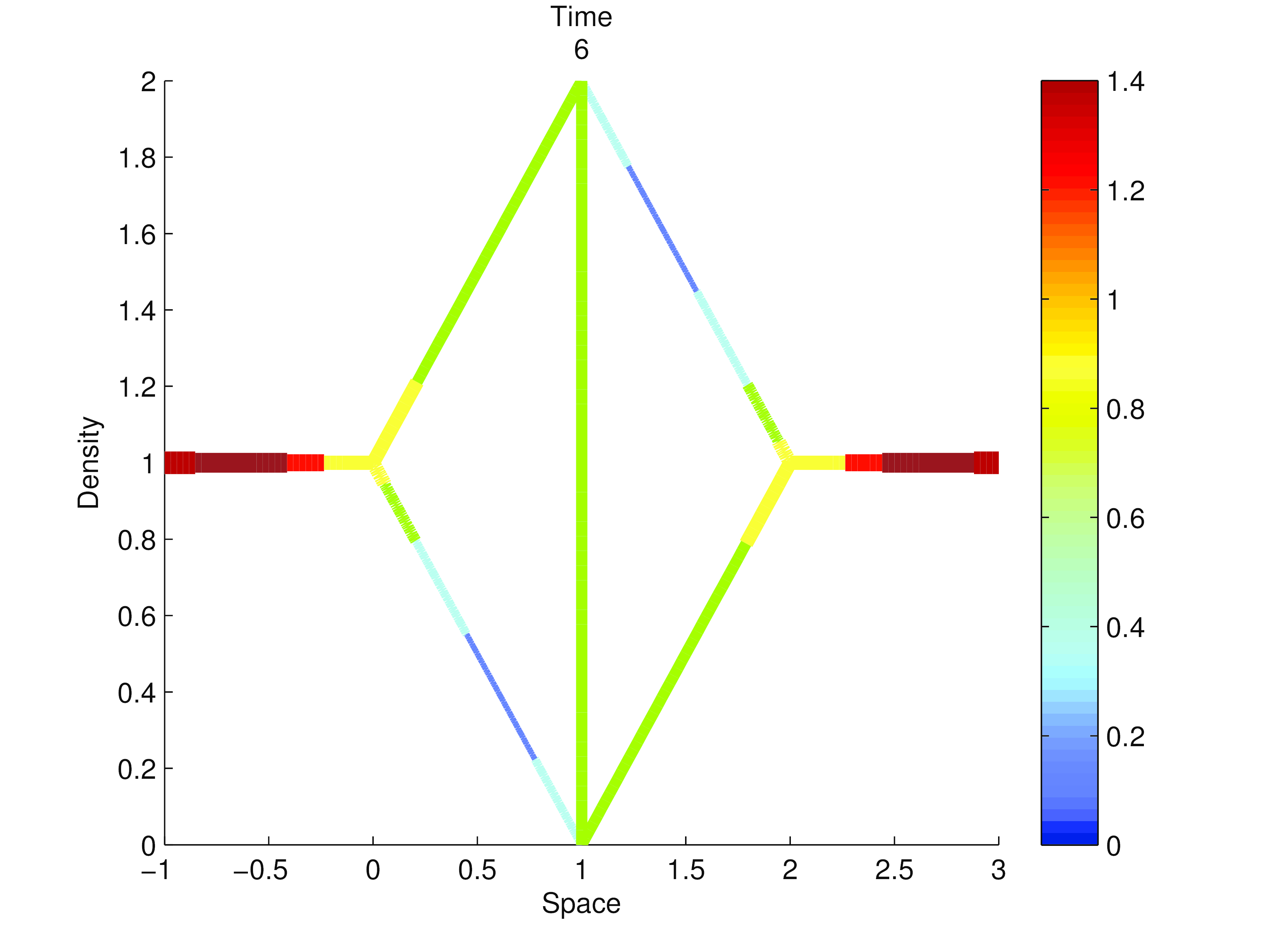}
\caption{Wheatstone bridge-shaped network. The distribution of the density $u_i(x,t)$ on each arc of the diamond graph-like network at time $t=2.5$ (on the left) and $t=6$ (on the right).}
\label{fig:2d}
\end{center}
\end{figure}

The effect on the movement of the organism is that it is able to find the shortest path connecting node $0$ with node $5$ in such a way to minimize the sum of the length of the arcs composing the path. As a result, see Figg. \ref{fig:1d}-\ref{fig:2d}, the mass concentration is higher on the arcs composing the shortest path, as underlined by the colourscale and thickness of the line connecting nodes.  
Note that the solution obtained at time $T=6$ is a stationary state. 

Then, to have a certain quantity of cells enters into the network by the node $0$ and the node $5$, we set non-zero flux conditions at the outer boundaries \eqref{inflow_cond_v1} at the node $0$ and \eqref{inflow_cond_v2} at the node $5$:
$$
v_1(0,t) = \frac{2}{1+ u_1(0,t)}, \quad v_{7}(L_{7},t) = -\frac{2}{1+ u_{7}(L_{7},t)}.
$$

\begin{figure}[htbp!]
\begin{center}
\includegraphics[height=5.5cm,width=5.5cm]{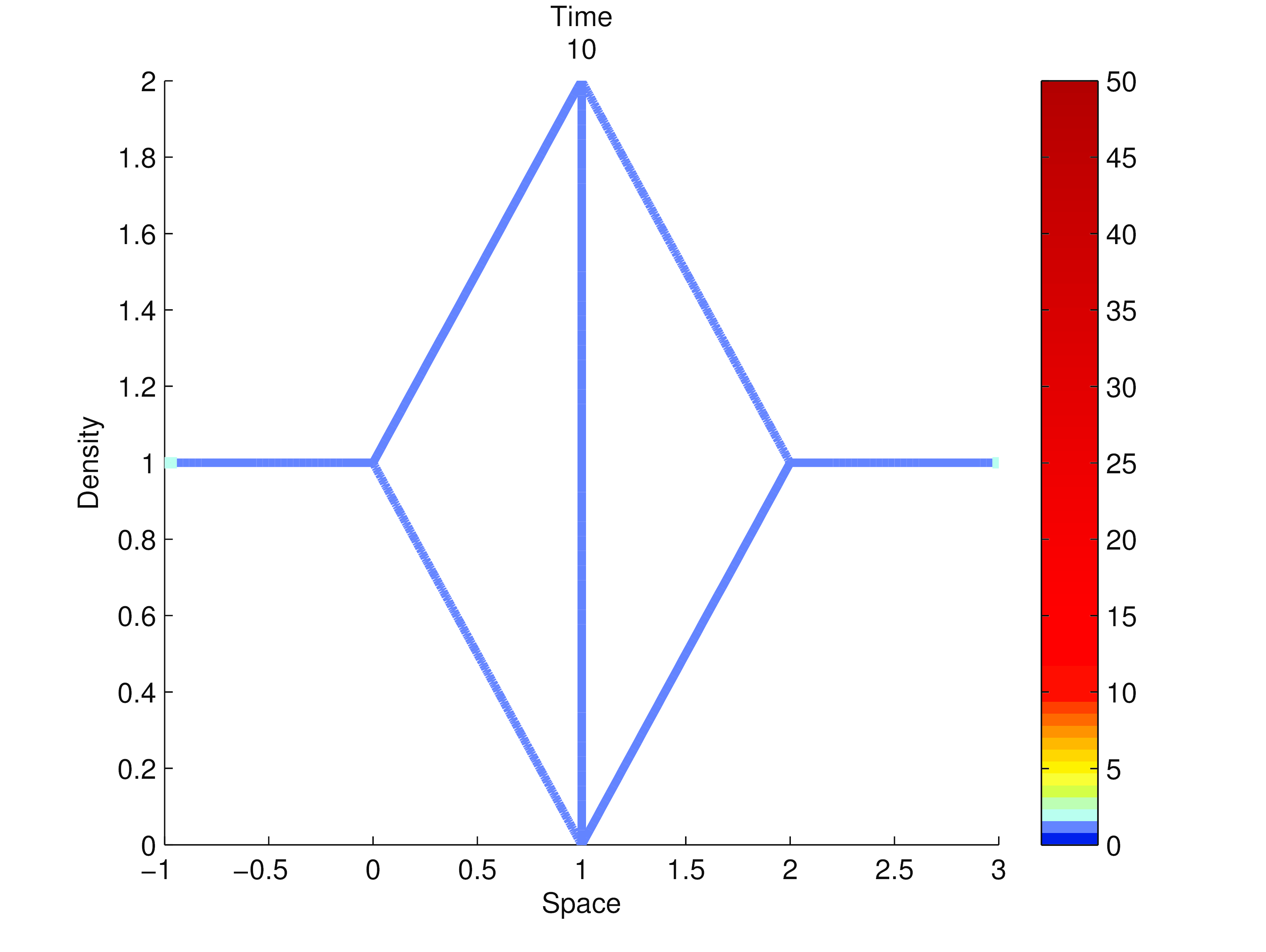}\quad \includegraphics[height=5.5cm,width=5.5cm]{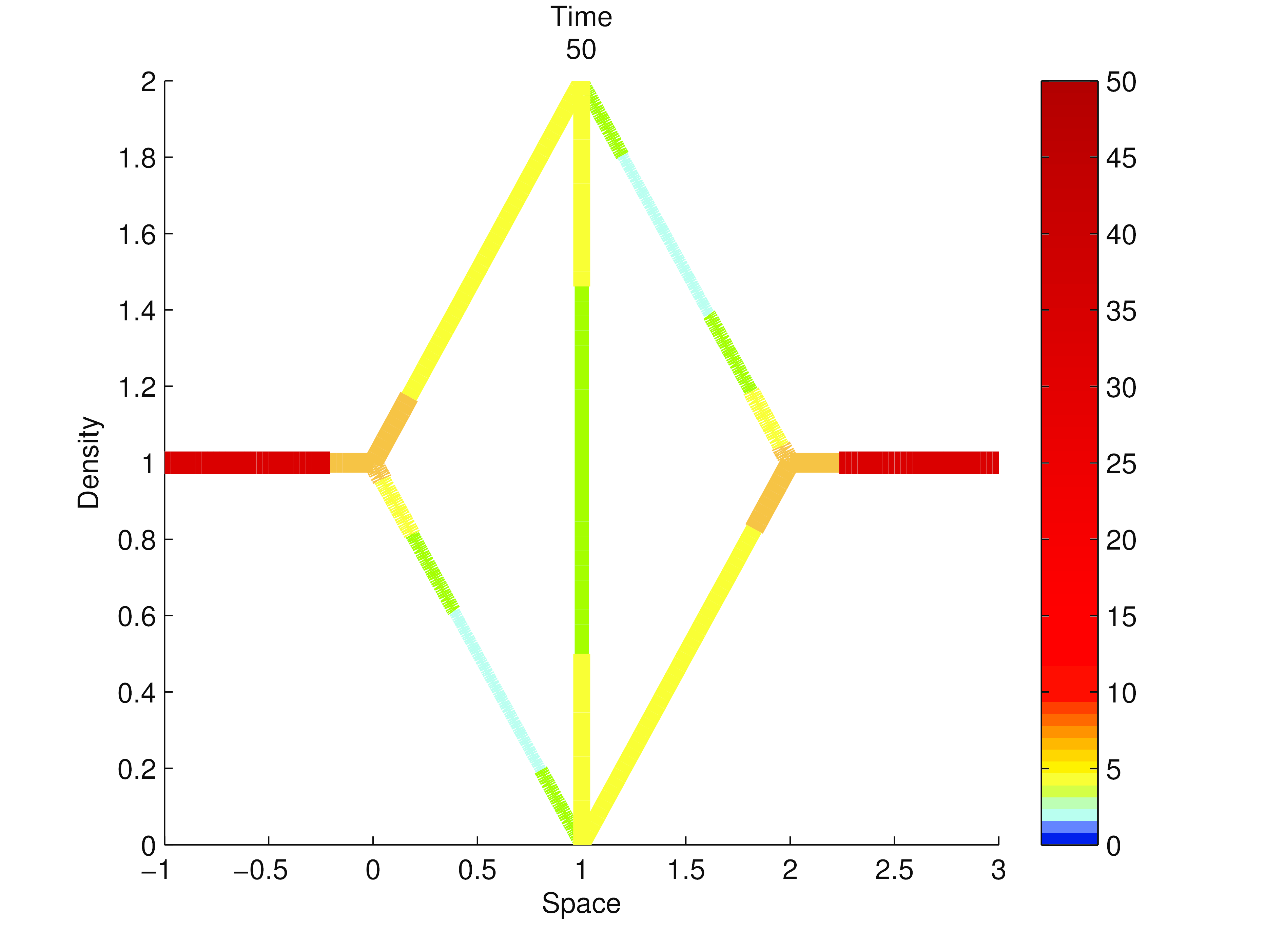}
\caption{Wheatstone bridge-shaped network. The distribution of the density $u_i(x,t)$ on each arc of the diamond graph-like network at time $10$ (on the left) and $t=50$ (on the right).}
\label{fig:3d}
\end{center}
\end{figure}

\begin{figure}[htbp!]
\begin{center}
\includegraphics[height=5.5cm,width=5.5cm]{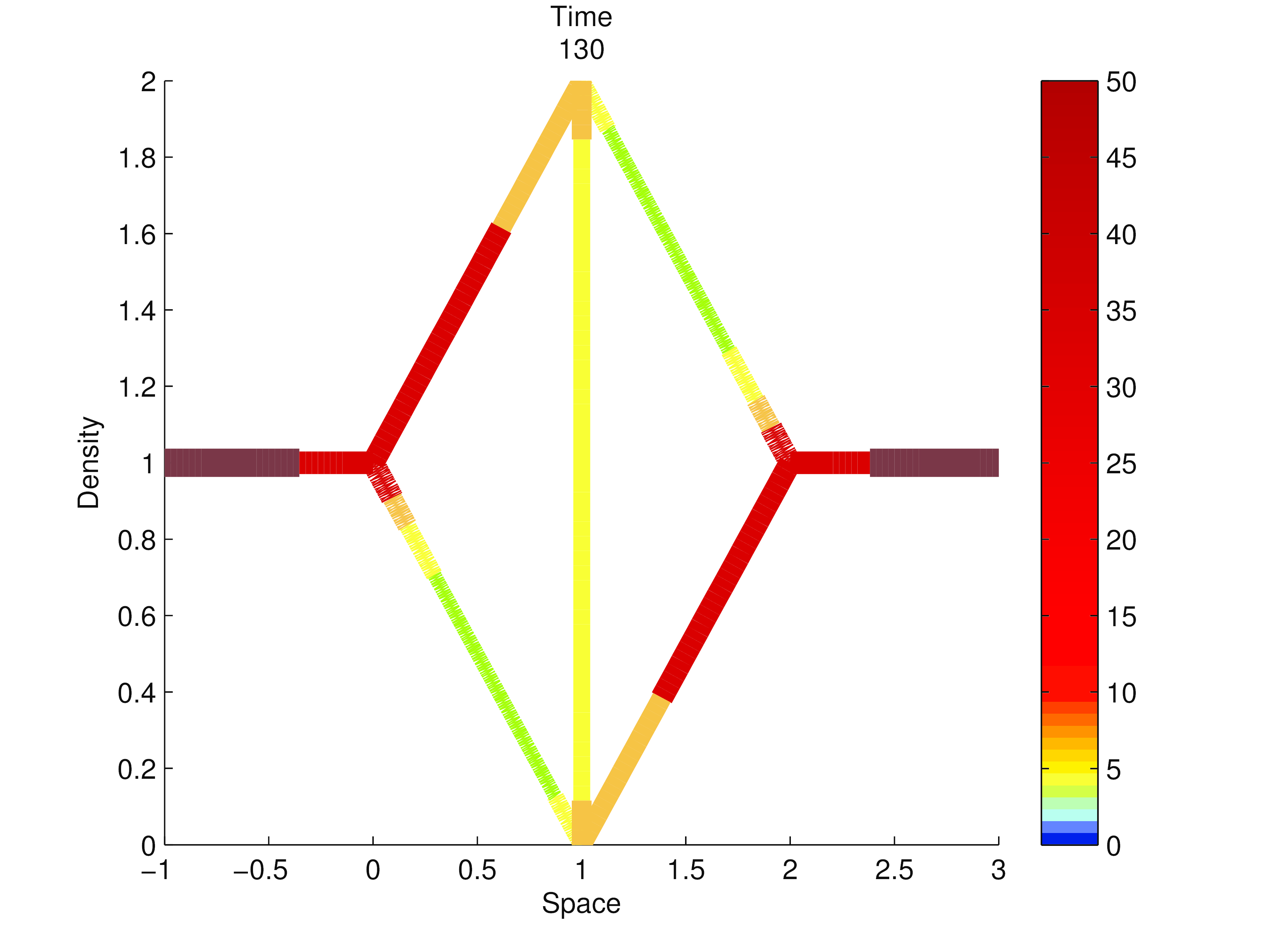}\quad \includegraphics[height=5.5cm,width=5.5cm]{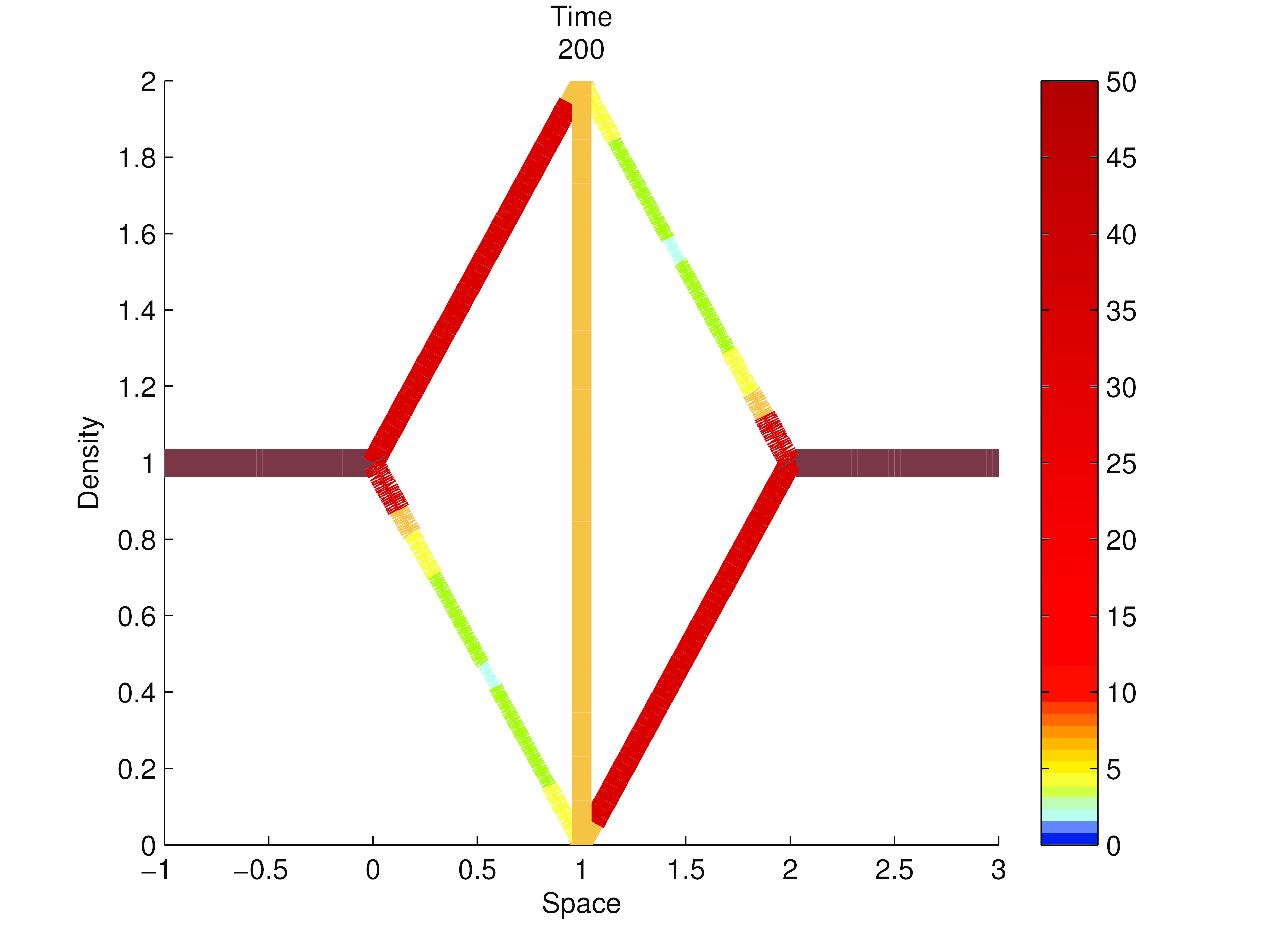}
\caption{Wheatstone bridge-shaped network. The distribution of the density $u_i(x,t)$ on each arc of the diamond graph-like network at time $t=130$ (on the left) and $t=200$ (on the right).}
\label{fig:4d}
\end{center}
\end{figure}

As in the previous case, cells migrate on the arcs on the shortest path connecting node $0$ with node $5$ in such a way to minimize the sum of the length of the arcs composing the path. As shown in the following Figures \ref{fig:3d}-\ref{fig:4d}, the mass concentration is higher on the path of minimum total length. In this case the connection between the two exits in more evident, as underlined by the range of the density in the colourscale.

\subsubsection{A network with two exits: the maze}

\begin{figure}[htbp!]
\begin{center}
\includegraphics[width=8cm]{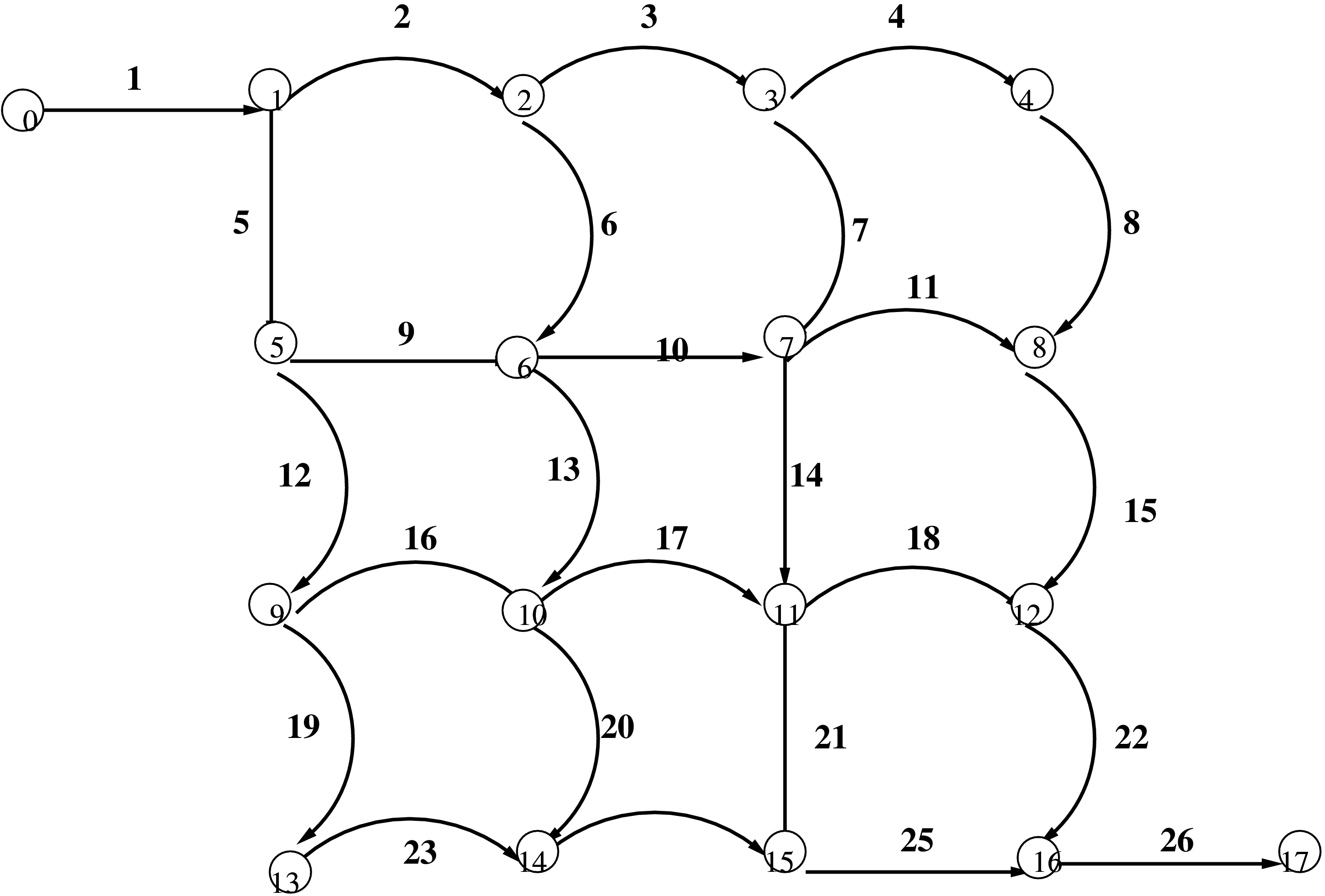}
\caption{The maze: a network with 26 arcs and 18 nodes.}
\label{fig:6}
\end{center}
\end{figure}

As in \cite{BGKS}, the network we consider here is a maze composed of 26 arcs and 18 nodes, with the exits placed at the node 0 (source) and at the node 17 (sink). We set the length $L_i=0.5$ on arcs $i=1,5,9,10,14,21,25,26$ and $L_i=10$ on the others.
Note that in Figure \ref{fig:6} we represent such network, with the shortest arcs depicted with the straight line and the longest arcs with the curve line.\\
  We set parameters as $a_i=1, b_i=0.1, \lambda_i=\sqrt{0.33}, D_i= 1, \chi_i= 1$, for all $i$, considering again positive chemotaxis. Furthermore, for the transmission conditions for $u$ we set dissipative coefficients $\xi_{i,j}$, see Table \ref{tab:coeff}, and for the transmission conditions on $\phi$ we assume again $\kappa_{i,j}=1$ for $i \neq j$ and $\kappa_{i,i}=0$, for $i,j=1,\ldots,26$.
The initial values for $u_i$ are randomly equally distributed in $[0.45,0.55]$ and we assume $\phi_i (x, 0) = 0$, $\forall i$. 

\begin{table}[h!]
\centering
{
\begin{tabular}{|c|c|} \hline
node p= 1, 2, 3, 5, 8, 9, 12, 14, 15, 16 &$\xi_{i,j}=\frac{1}{3}, i,j \in I_p \cup  O_p$ \\  \hline

node p= 4, 13 &$\xi_{i,j}=\frac{1}{2},  i,j \in I_p \cup  O_p$\\  \hline

node p= 6, 7, 10, 11 &$\xi_{i,j}=\frac{1}{4},  i,j \in I_p \cup  O_p$ \\  \hline
\end{tabular}}
\vspace{0.1 in} 
\caption{Transmission coefficients used for the numerical simulations of the maze in Fig. \ref{fig:6} given node by node.}
\label{tab:coeff}
\end{table}

Here we set non-homogeneous Neumann boundary conditions at the outer boundaries for $\phi$: inflow conditions leading to a high concentration of the chemoattractant at the specified point are assumed at the node $0$ and at the node $17$:

\begin{equation}\label{inflow_cond_phi:1} 
\partial_x \phi_1(0,t) = -1 , \quad \partial_x \phi_{26}(L_{26},t).
\end{equation}

Furthermore, to have a certain quantity of cells enters into the network by the node $0$ and the node $17$, we set non-zero flux conditions at the outer boundaries:
$$
v_1(0,t) = \frac{2}{1+ u_1(0,t)}, \quad v_{26}(L_{26},t) = -\frac{2}{1+ u_{26}(L_{26},t)}.
$$
The effect on the movement of the cells is that they try to find the  minimum-length way connecting node $0$ with node $17$ in order to minimize the sum of the length of the arcs composing the path. Let us recall that the dead-end cutting of the zones with no-flux and the shortest path selection is here represented, as shown in the following Figures \ref{fig:8}-\ref{fig:9}, by the higher mass concentration on the edges composing the shortest path.  
We also observe that at a certain time blow up of solutions verifies, due to the growth of the total mass. We remark that, though the evolution is slower respect to that presented in \cite{BGKS}, the results are asymptotically the same obtained in \cite{BGKS}.

\begin{figure}[htbp!]
\begin{center}
\includegraphics[height=5.5cm,width=5.5cm]{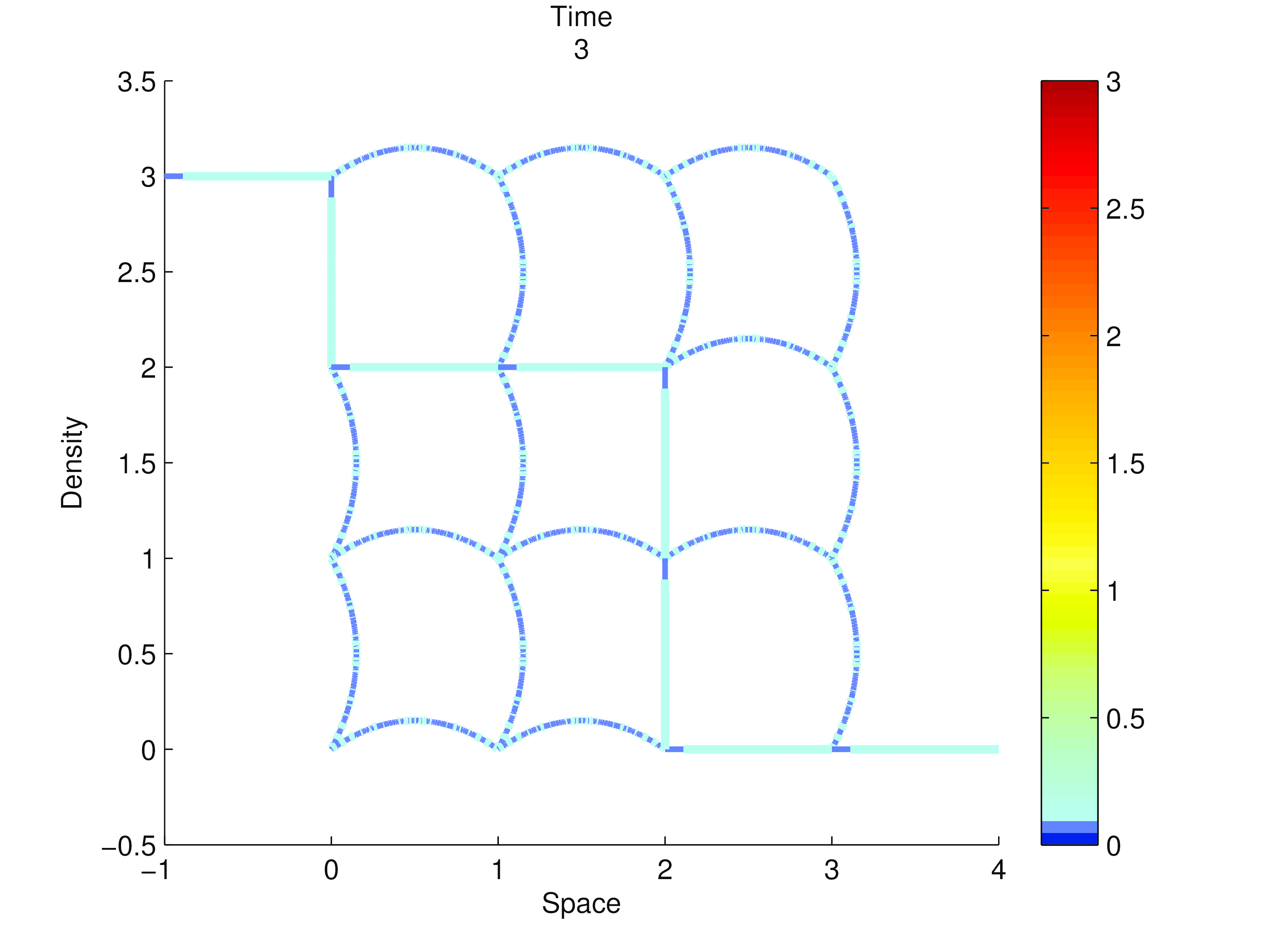}\quad \includegraphics[height=5.5cm,width=5.5cm]{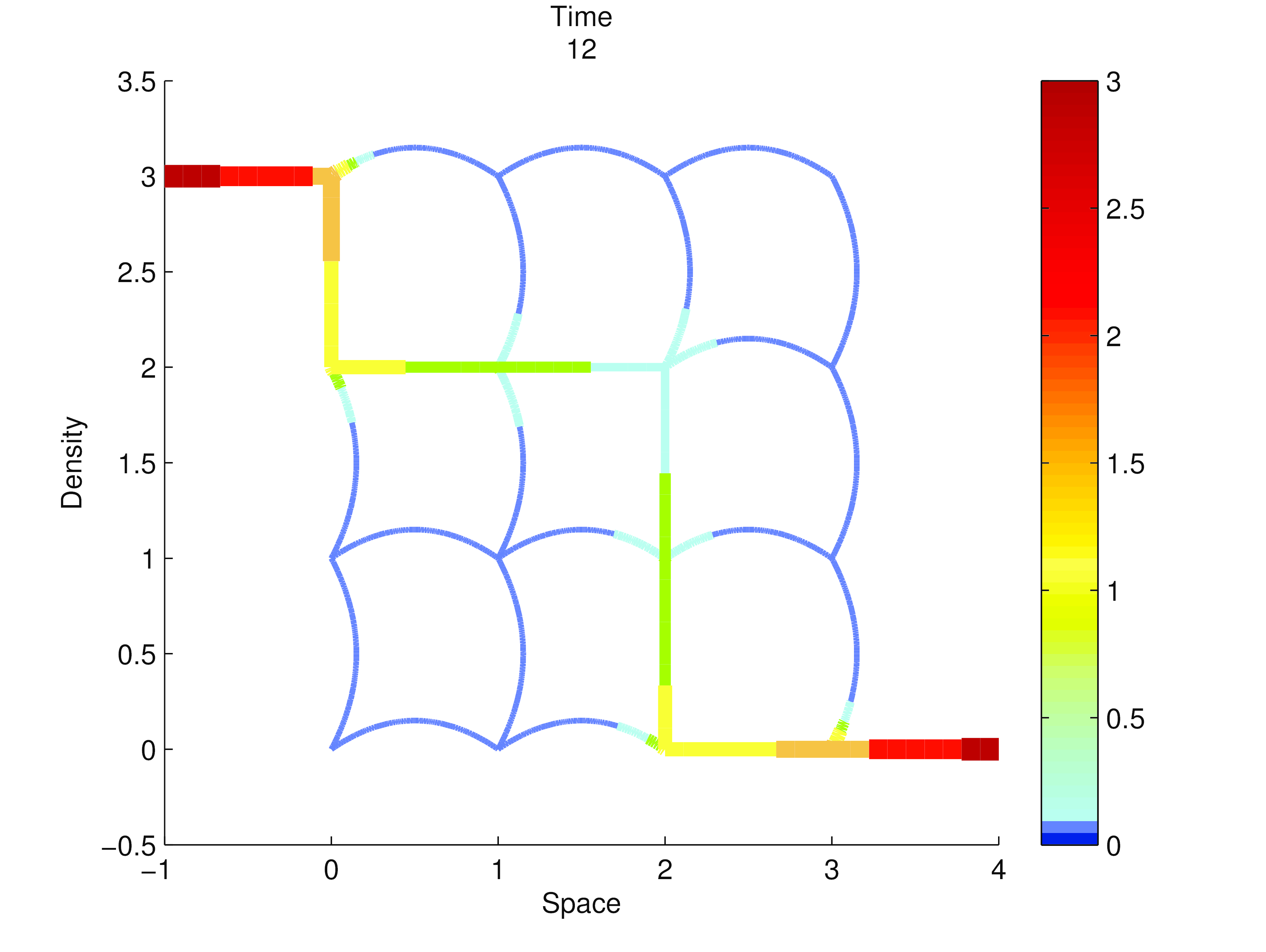}
\caption{The maze. The distribution of the density $u_i(x,t)$ on each arc of the network in Fig. \ref{fig:6} at time $t=3$ (on the left) and $t=12$ (on the right).}
\label{fig:8}
\end{center}
\end{figure}
\begin{figure}[htbp!]
\begin{center}
\includegraphics[height=5.5cm,width=5.5cm]{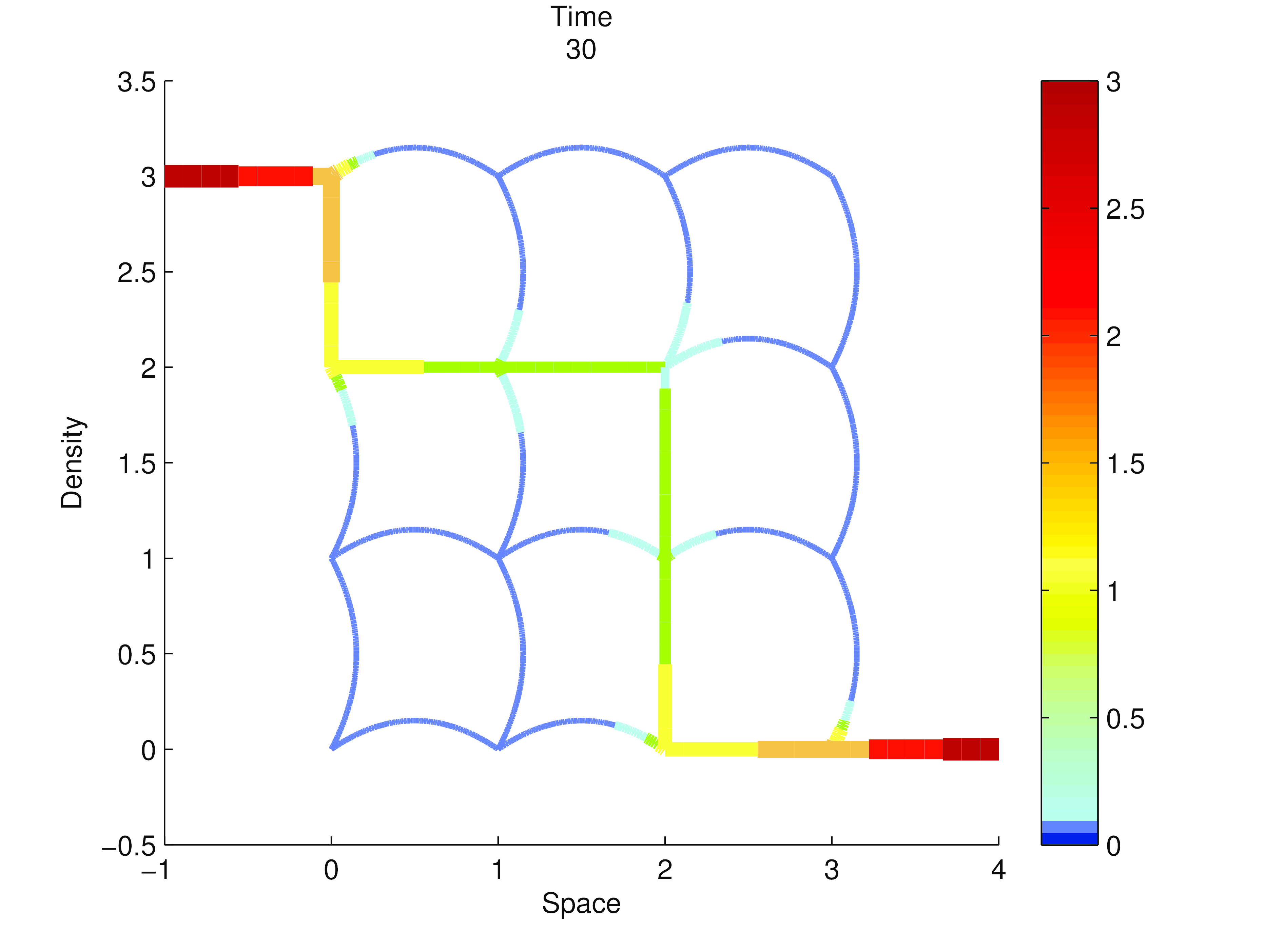}\quad \includegraphics[height=5.5cm,width=5.5cm]{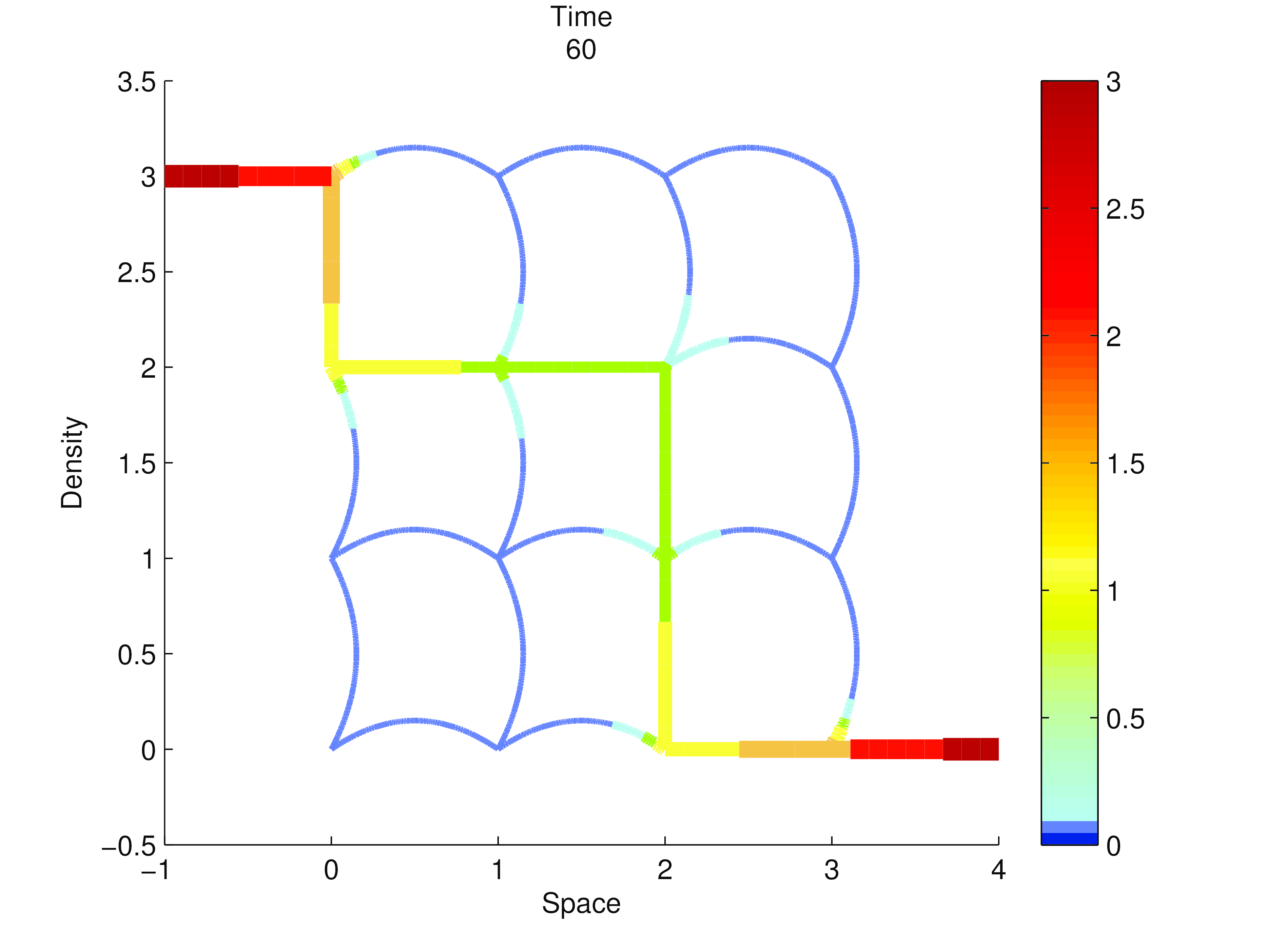}
\caption{The maze. The distribution of the density $u_i(x,t)$ on each arc of the network in Fig. \ref{fig:6} at time $t=30$ (on the left) and $t=60$ (on the right).}
\label{fig:9}
\end{center}
\end{figure}
The simulation of the network above was performed by a personal computer, processor Intel Centrino core 2 Duo 2 Ghz, RAM 3 Gb and the CPU time for 1333 time iterations was 19.71 s.

\subsubsection{A network with more than two exits}
Here we consider a network with more than two exits (NM2E): in fact in this case we have two sources and three sinks.
Such network, which  is composed of 21 arcs and 15 nodes, with the sources placed at the nodes 0 and 1, while the sinks are at the nodes 9, 10 and 11, we set the length $L_i=0.5$ on arcs $i=1,2,3,4,8,13,14,15,19,20,21$ and $L_i=10$ on the others. The network we consider is shown in Figure \ref{fig:6}, where the shortest arcs are depicted with the dashed line.\\
As before, we set parameters as $a_i=1, b_i=0.1, \lambda_i=\sqrt{0.33}, D_i= 1, \chi_i= 1$, for all $i$. Furthermore, for the transmission conditions for $u$, we set dissipative coefficients $\xi_{i,j}$, see Table \ref{tab:coeff2}, and for the transmission conditions on $\phi$ we assume $\kappa_{i,j}=1$ for $i \neq j$ and $\kappa_{i,i}=0$, for $i,j=1,\ldots,21$.
We assume the inflow boundary conditions:
\begin{equation}\label{inflow_cond_phi} 
\partial_x \phi_1(0,t) = -1 = \partial_x \phi_2(0,t) , \quad \partial_x \phi_{19}(L_{19},t) = 1 = \partial_x \phi_{20}(L_{20},t) = \partial_x \phi_{21}(L_{21},t).
\end{equation}
Moreover, we assume for $u$ homogeneous Neumann boundary conditions $u_x(\cdot,t)=0$ at the nodes $0, 1$ and at the nodes $9, 10, 11$.

\begin{figure}[htbp!]
\begin{center}
\includegraphics[width=8cm]{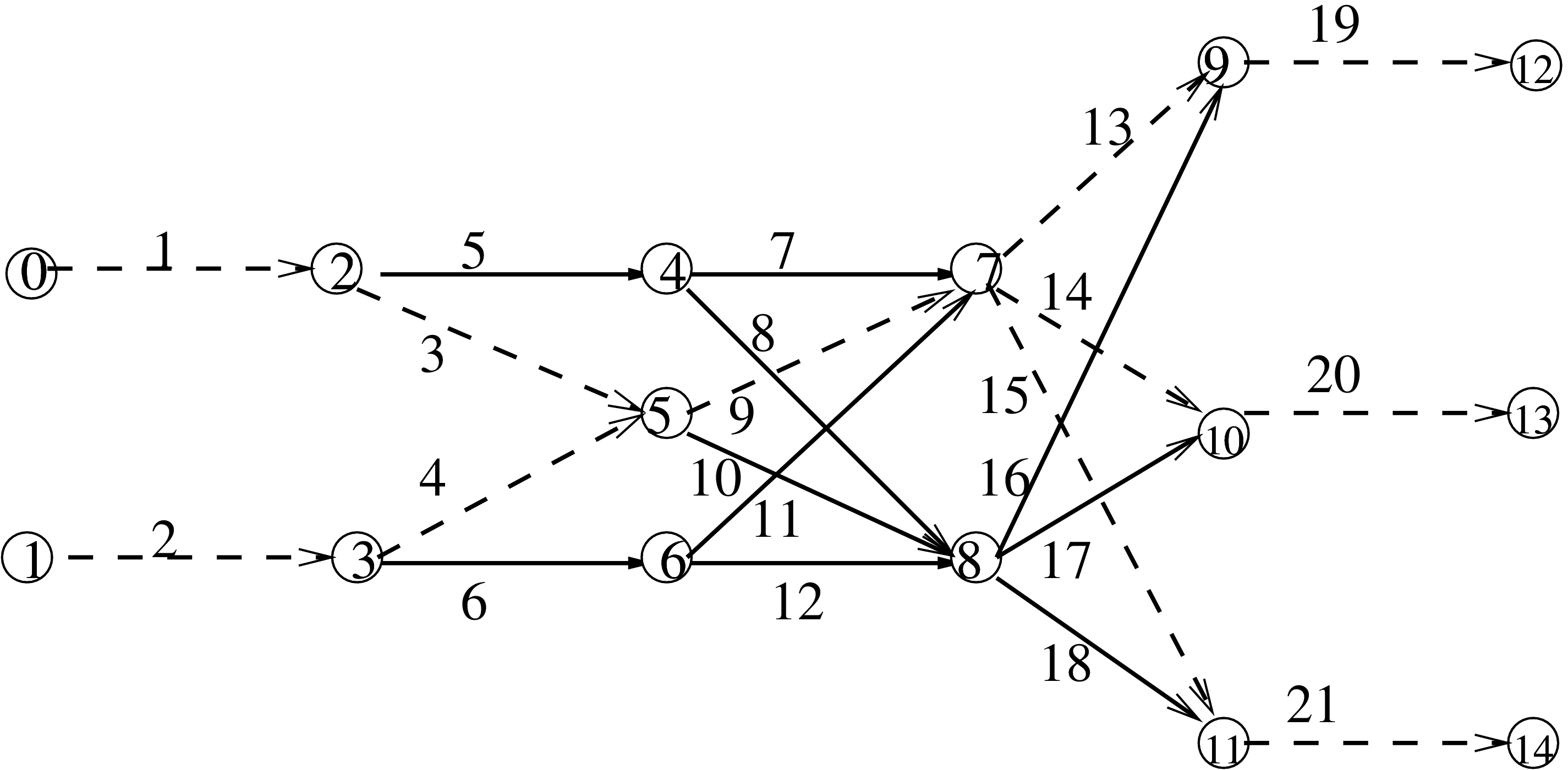}
\caption{Network of 21 arcs and 15 nodes (NM2E): the dashed line represents the shortest arcs.}
\label{fig:10}
\end{center}
\end{figure}
The initial values for $u_i$ are randomly equally distributed in $[0.45,0.55]$ and we assume $\phi_i (x, 0) = 0$, $\forall i$.

\begin{table}[h!]
\centering
{
\begin{tabular}{|c|c|} \hline
node p=  2, 3, 4, 6, 9, 10, 11 &$\xi_{i,j}=\frac{1}{3}, i,j \in I_p \cup  O_p$ \\  \hline

node p= 7, 8 &$\xi_{i,j}=\frac{1}{6},  i,j \in I_p \cup  O_p$\\  \hline

node p= 5 &$\xi_{i,j}=\frac{1}{4},  i,j \in I_p \cup  O_p$ \\  \hline
\end{tabular}}
\vspace{0.1 in} 
\caption{Transmission coefficients used for the numerical simulations of network NM2E in Fig. \ref{fig:6} given node by node.}
\label{tab:coeff2}
\end{table}
\begin{figure}[htbp!]
\begin{center}
\includegraphics[height=5.5cm,width=5.5cm]{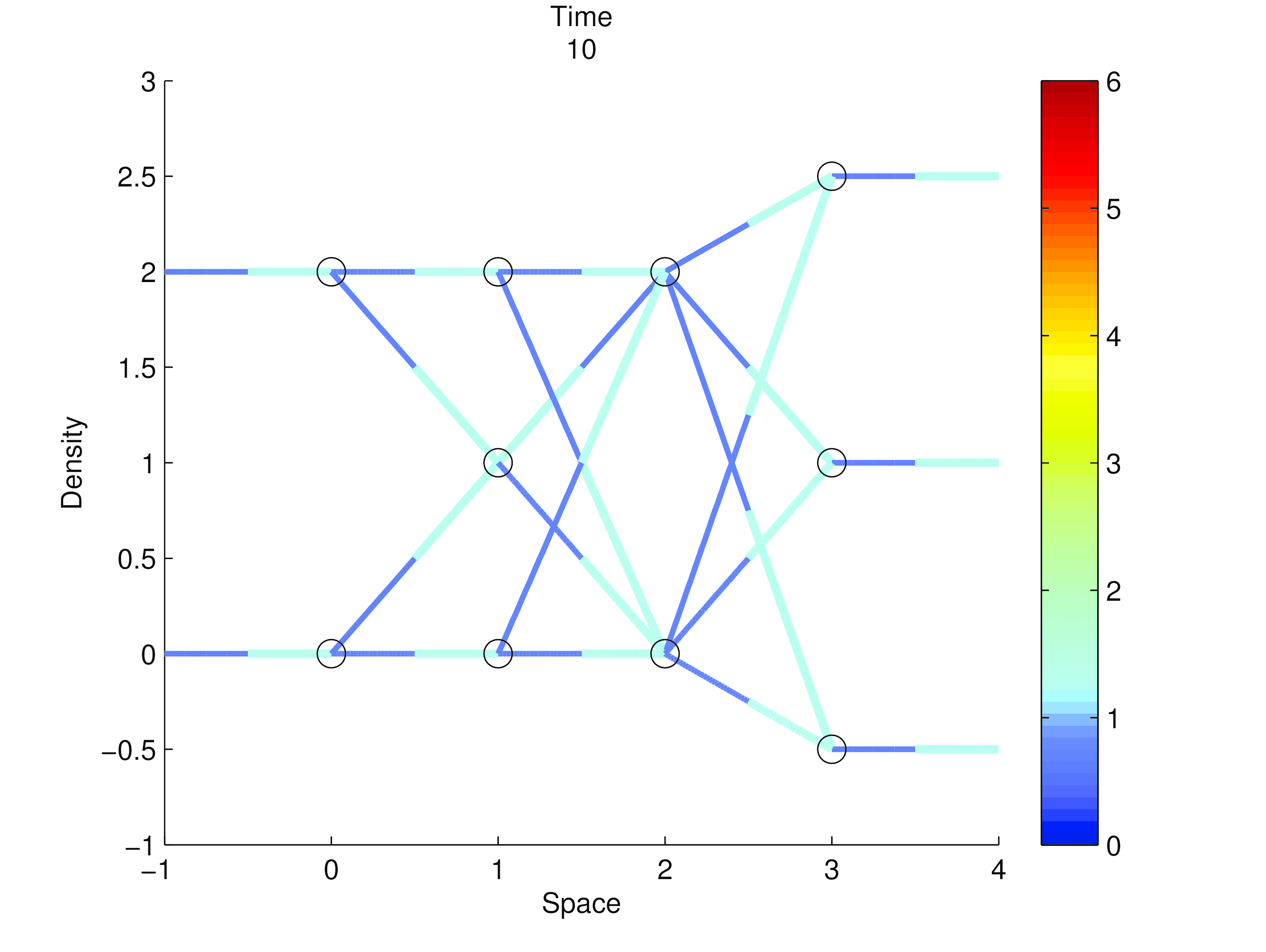}\quad \includegraphics[height=5.5cm,width=5.5cm]{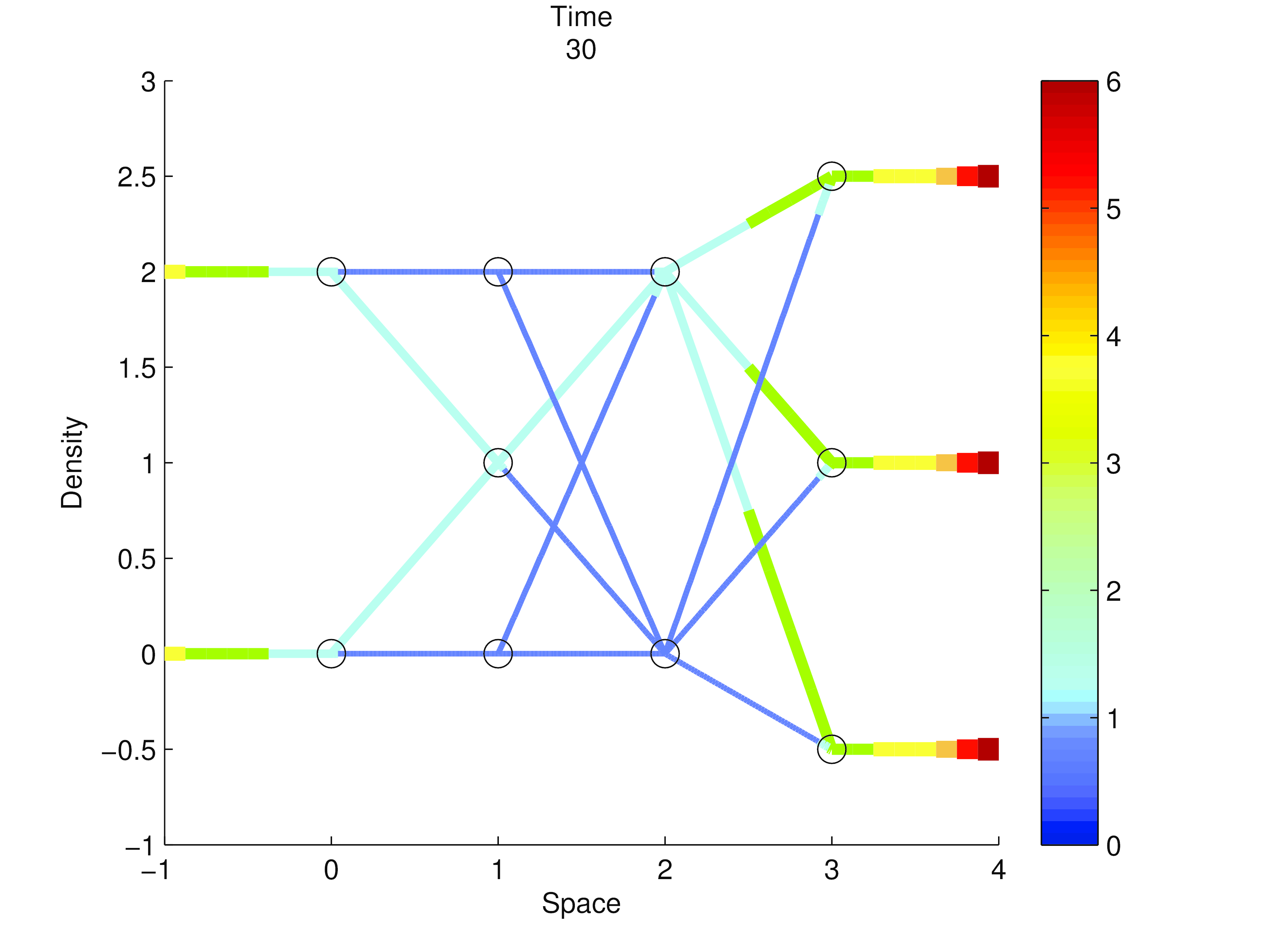}
\caption{The distribution of the density $u_i(x,t)$ on each arc of the network NM2E in Fig. \ref{fig:10} at time $t=10$ (on the left) and $t=30$ (on the right).}
\label{fig:11}
\end{center}
\end{figure}
\begin{figure}[htbp!]
\begin{center}
\includegraphics[height=5.5cm,width=5.5cm]{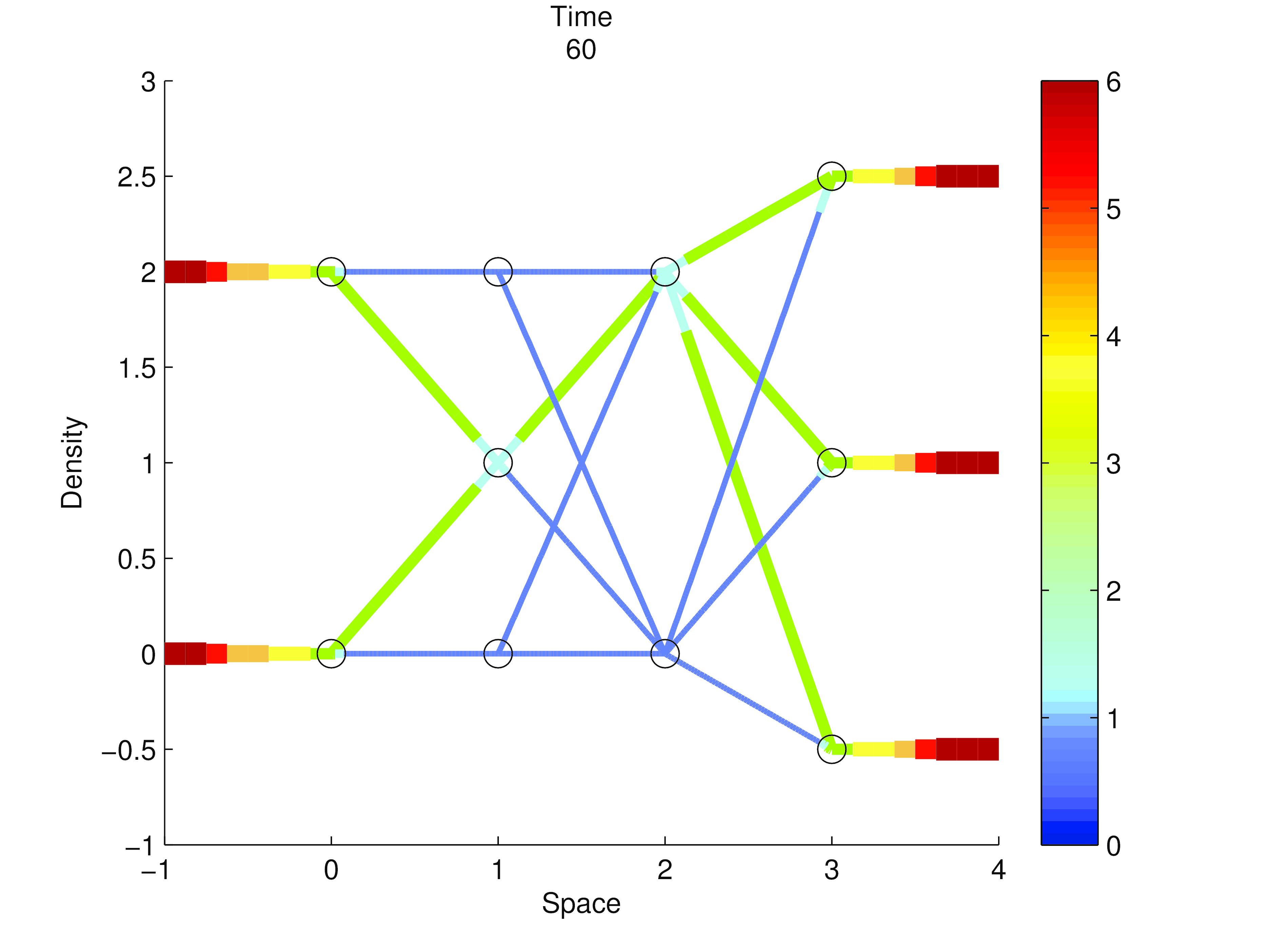}\quad \includegraphics[height=5.5cm,width=5.5cm]{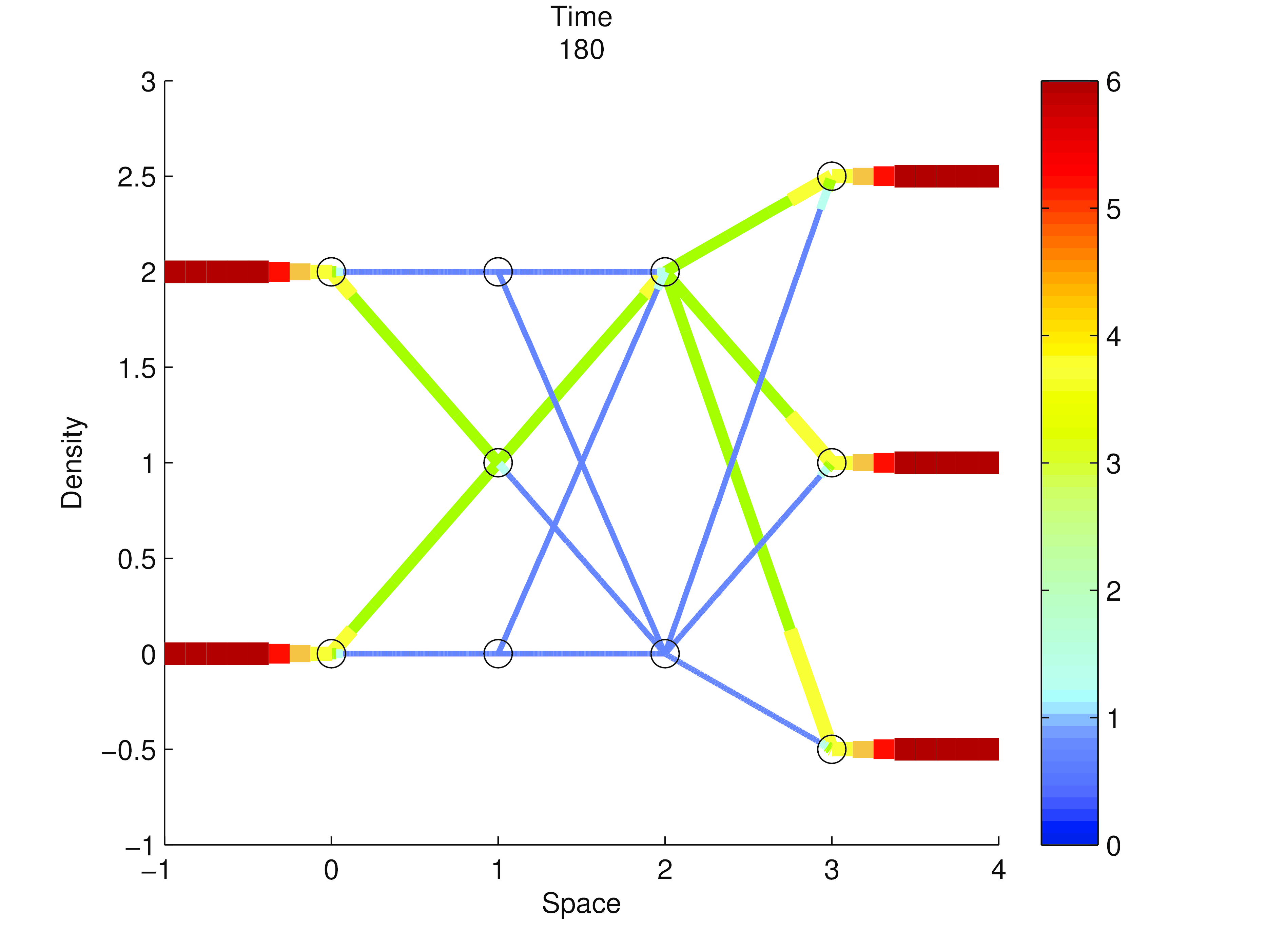}
\caption{The distribution of the density $u_i(x,t)$ on each arc of the network NM2E in Fig. \ref{fig:10} at time $t=60$ (on the left) and $t=180$ (on the right).}
\label{fig:12}
\end{center}
\end{figure}

Even in this case, the organism traces the path connecting nodes $0$ and $1$ with nodes $19, 20, 21$ in such a way to minimize the sum of the total length of the arcs. As a result, see the following Figures \ref{fig:11}-\ref{fig:12}, besides the arcs  connected to the exits ($1,2,19,20,21$), the path includes the arcs  with minimum length $3,4,8,13,14,15$.

\end{document}